\documentclass[3p]{elsarticle}
\usepackage[utf8]{inputenc}
\usepackage{graphicx}
\usepackage{amsmath}
\usepackage{amssymb}
\usepackage{amsfonts}
\usepackage{float} 
\usepackage{pdfpages} 
\usepackage{subfigure}
\usepackage{wrapfig} 
\usepackage{textcomp} 
\usepackage{multirow} 
\usepackage{microtype}
\usepackage{tabularx} 
\usepackage{rotating}
\usepackage{natbib}
\usepackage[footnotesize]{caption}
\usepackage{hhline}
\usepackage{supertabular}
\usepackage{enumitem}
\usepackage{array}
\usepackage{color}
\usepackage{caption}

\usepackage{hyperref}
\hypersetup{
    colorlinks,%
    citecolor=black,%
    filecolor=black,%
    linkcolor=black,%
    urlcolor=black}

\DeclareGraphicsExtensions{.pdf, .png, .jpg}
\graphicspath{{graphs/}} 

\journal{Applied Mathematics and Computation}

\begin{document}

\begin{frontmatter}


\newtheorem{mydef}{Definition}
\newtheorem{proposition}{Proposition}
\newtheorem{remark}{Remark}
%

\title{A numerical method for junctions in networks of shallow-water channels}

\author[unitn]{Francesca Bellamoli}
\author[ntnu]{Lucas O. M\"{u}ller\corref{cor1}}
\ead{lucas.o.muller@ntnu.no}
\cortext[cor1]{Corresponding author}

\author[unitn]{Eleuterio F. Toro}

\address[unitn]{Laboratory of Applied Mathematics,DICAM,\\
        University of Trento,\\
        Trento, Italy.}
\address[ntnu]{Division of Biomechanics, Department of Structural Engineering,\\
	  Norwegian University of Science and Technology,\\
	  Trondheim, Norway.}

\begin{abstract}

There is growing interest in developing mathematical models and appropriate numerical methods for problems involving networks formed by, essentially, one-dimensional (1D) domains joined by junctions.  Examples include hyperbolic equations in networks of gas tubes, water channels and vessel networks for blood and lymph in the human circulatory system. A key point in designing numerical methods for such applications is the treatment of junctions, \emph{i.e.} points at which two or more 1D domains converge and where the flow exhibits multidimensional behaviour. This paper focuses on the design of methods for networks of water channels. Our methods adopt the finite volume approach to make full use of the two-dimensional shallow water equations on the true physical domain, locally at junctions, while solving the usual one-dimensional shallow water equations away from the junctions. In addition to mass conservation, our methods enforce conservation of momentum at junctions;  the latter seems to be the missing element in methods currently available. Apart from simplicity and robustness, the salient feature of the proposed methods is their ability to successfully deal with transcritical and supercritical flows at junctions, a property not enjoyed by existing published methodologies. Systematic assessment of the proposed methods for a variety of flow configurations is carried out. The methods are directly applicable to other systems, provided the multidimensional versions of the 1D equations are available.

\end{abstract}


\end{frontmatter}

\section{Introduction}

There are multidimensional physical problems modelled by partial differential equations in networks of spatial domains than are essentially straight. In such cases the governing equations can be assumed to be one-dimensional (1D), potentially resulting in significant computing savings. Examples include gas flow in pipes 
\cite{Banda:2006a, Brouwer:2011a, Bales:2009a,Reigstad:2015a,Bermudez:2017},
traffic flow \cite{Coclite:2005a, Borsche:2014c, Bretti:2007a} 
water flows \cite{Kesserwani,AkanYen,AralZhangJin,Zhang, Borsche:2014b, Kesserwani:2008a} 
and blood flow in the human circulatory system 
\cite{Quarteroni:2000a,Formaggia:1999a,Olufsen:2000a, Sherwin2003,Miglio2005_1,Alastruey:2011a,Liang:2009a,Liang:2009b,FullanaZaleski,Liang:2014a,Mueller:2014a,Mueller:2014b,Toro:2015a}. The challenge, however, is how to connect these 1D domains in a way that accounts for the multidimensional character of the equations, even in an approximately manner. 

Current methods are reported to perform well in most cases. However, a shortcoming of existing methods is their inability to deal with transcritical and supercritical flow through junctions. In some cases, these methods fail even for subcritical flows at moderately high Froude numbers. Transcritical and supercritical flows are important flow regimes that may occur more often than one is aware of, for example at junctions, locally. In physiological flows this is found to be the case in the venous system, under postural changes. In open channel flows the occurrence of supercritical flows is not rare and may potentially take place in inundating flows emerging from dam-break events and tsunami waves. Supercritical regimes may also appear in networks of tubes transporting compressible fluids under extreme accidental conditions. 

In this paper we present methods for dealing with junctions connecting 1D domains and illustrate the ideas for junctions of 1D shallow water channels. We note that the full problem is governed by the two-dimensional (2D) shallow water equations. The methods presented here make use of the finite volume approach, whereby the true geometry is accounted for locally at junctions, whereas away from junctions, the usual 1D equations are solved. In addition to mass conservation, our methods enforce conservation of momentum at junctions, which constitutes an improvement over methods currently available. It is noted that the approach, as applied to complex networks of channels, can lead to very significant computing savings, as compared to solving the full multidimensional problem, without compromising the solution quality.  Systematic assessment of the methods for a variety of flow configurations is carried out. 
It is worth noting that a similar method, which combines a 1D model for the channels and a 2D model on unstructured grids for junctions, has been investigated both in hydrodynamics (Miglio et al. \cite{Miglio2005_1}) and in haemodynamics  (Formaggia et al \cite{Formaggia1999,FormaggiaQuarteroni2003}). However, Miglio et al. (\cite{Miglio2005_1}) use a finite element scheme and investigate only the case of subcritical flows.

The rest of this paper is structured as follows. In Section \ref{sec:mathmodnumerics} we briefly present the underlying mathematical models and the numerical framework upon which the proposed methods are constructed. Next, the novel methodology for coupling 1D domains at junctions is illustrated in Section \ref{sec:method}. Section \ref{sec:results} is devoted to the validation of the proposed methods, while  concluding remarks are presented in Section \ref{sec:discussion}. \ref{section_1Dexisting} describes an existing method for junctions and \ref{chapter_theoretical} provides details on the 2D unstructured mesh method used to produce reliable reference solutions used to assess the methods presented in this paper. 

\section{Mathematical models and numerical method} \label{sec:mathmodnumerics}

\begin{figure}[H]
\centering
\includegraphics[width=0.9\textwidth]{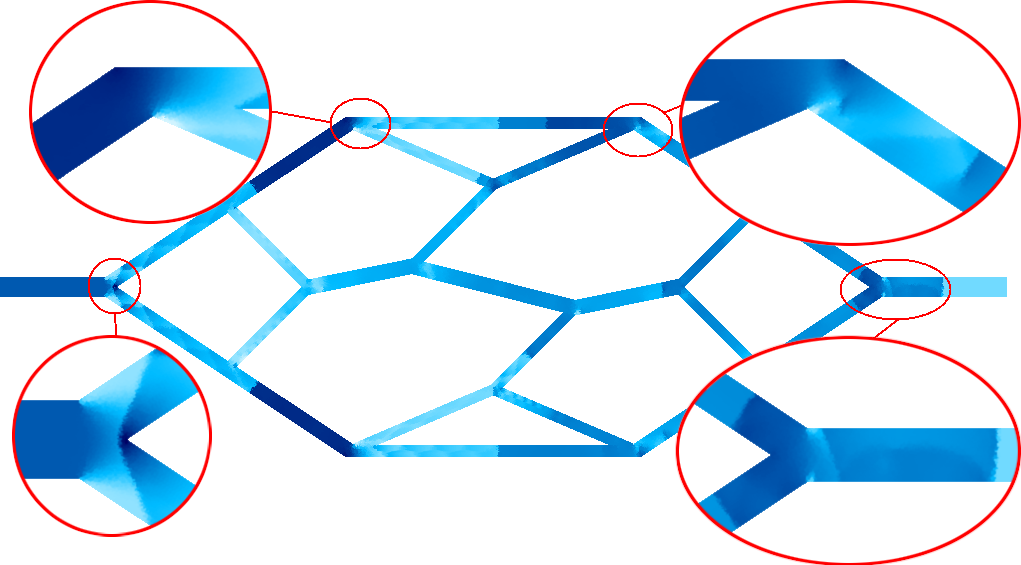}
\caption{Example of a channel network. Regions with two-dimensional behaviour are encircled and zoomed in. }
\label{introduzione_network}
\end{figure}

We are concerned with free-surface shallow water flow under gravity in a network system consisting of interconnected straight (or essentially straight) channels joined at junctions, as illustrated in Fig. \ref{introduzione_network}. The flow field has a significant 2D structure only in the vicinity of junctions, while it is essentially 1D along the straight channels, away from junctions. The purpose of this work is to develop a method that combines the use of a 1D model in channels and a 2D model only at junctions, coupling these two models with appropriate matching conditions. As we shall show later, the resulting methods show a huge computational efficiency gain with respect to solving the full 2D equations on the entire domain. 
The methods proposed in this paper adopt the finite volume framework and combine the 2D equations and its corresponding local geometry in a neighbourhood of junctions, along with the 1D equations in the straight channels. Next, we recall the governing equations and the finite volume method.

\subsection{The shallow water equations}

The time-dependent, non-linear, 2D shallow water equations written in conservation form read
\begin{equation} \label{swe1}
\partial_{t}{\bf Q} +  \partial_{x}{\bf F}({\bf Q})  +  \partial_{y}{\bf G}({\bf Q}) = {\bf S} ({\bf Q}) \;,
\end{equation}
with
\begin{equation} \label{swe2}
\left.
\begin{array}{ccc}
      {\bf Q}  
      = \left[ \begin{array}{c}
      h \\
      h u  \\
      h v   \\
      \end{array} \right] \;,\hspace{1mm} &
      {\bf F}({\bf Q}) 
      = \left[ \begin{array}{c}
      h u  \\
      h u^{2}+\frac{1}{2} g h^{2} \\
      h u v 
      \end{array} \right]  \;,  \hspace{2mm}  &
      {\bf G}({\bf Q}) 
      = \left[ \begin{array}{c}
      h v  \\
      h v u  \\
      h v^{2}+\frac{1}{2} g h^{2}
      \end{array} \right]   \;, \\
      \\
      & {\bf S}({\bf Q})
      = \left[ \begin{array}{c}
      0  \\
       gh(S_{ox}-S_{fx})  \\
       gh(S_{oy}-S_{fy})  \\
      \end{array} \right] \:. &
\end{array} \right\}
\end{equation} 
Here ${\bf Q}$ is the vector of conserved variables; ${\bf F}({\bf Q})$ and ${\bf G}({\bf Q})$ are the fluxes in the $x$ and $y$ directions, respectively, and ${\bf S} ({\bf Q})$ is the vector of source terms. The physical variables are water depth $h(x,y,t)$ and velocity components $u(x,y,t)$ and $v(x,y,t)$, in the $x$ and $y$ directions respectively. In this paper the source term vector ${\bf S}({\bf Q})$ accounts for the variation of the bottom topography
\begin{equation} \label{swe3}
\begin{array}{ccc}
    S_{ox}=-\partial_xb(x,y) & \quad\mbox{and}\quad\; & S_{oy}=-\partial_yb(x,y)
\end{array}
\end{equation} 
and the bed friction
\begin{equation}\label{swe4}
\begin{array}{ccc}
     S_{fx}=\displaystyle{\frac{n^2\,u\,\sqrt{u^2+v^2}}{h^{4/3}}} & \quad\mbox{and}\quad \;
     & S_{fy}=\displaystyle{\frac{n^2\,v\,\sqrt{u^2+v^2}}{h^{4/3}}}\;.
\end{array}
\end{equation} 
Here $b(x,y)$ represents bottom elevation above a horizontal datum, $n$ is the Manning's coefficient and $g$ is the acceleration due to gravity. In this work we only consider channels with horizontal bottom, that is $S_{ox}=0$, $S_{oy}=0$. Equations (\ref{swe1}) form a system of partial differential equations of hyperbolic type. For background on the shallow water equations and associated numerical methods see, for instance, \cite{Toro2001} and the many references therein.

The one-dimensional version of (\ref{swe1}), in the generic $s$ direction, reads
\begin{equation}\label{swe5}
     \partial_{t}{\bf Q}_{1D} +  \partial_{s}{\bf F}_{1D}({\bf Q}_{1D}) = {\bf S}_{1D}({\bf Q}_{1D})  \;,
\end{equation}
with
\begin{equation}\label{swe6}
\begin{array}{lll}
      {\bf Q}_{1D} 
      = \left[ \begin{array}{c}
      h_{1D} \\
      h_{1D} u_{1D}  \\
      \end{array} \right] \;,\hspace{3mm} &
      {\bf F}_{1D}({\bf Q}_{1D}) 
      = \left[ \begin{array}{c}
      h_{1D} u_{1D}  \\
      h_{1D} u^{2}_{1D}+\frac{1}{2} g h^{2}_{1D} \\
      \end{array} \right] \;,\hspace{3mm} &
        {\bf S}_{1D}({\bf Q}_{1D})
      = \left[ \begin{array}{c}
      0  \\
       gh_{1D}(S_{os}-S_{fs})  \\
      \end{array} \right]  \;,
\end{array}
\end{equation} 
where $h_{1D} = h$ and $u_{1D}$ is the velocity along $s$.
\subsection{Rotational invariance and numerical method}

\begin{figure}
       \centerline{
       \includegraphics[scale=0.40,angle=-90]{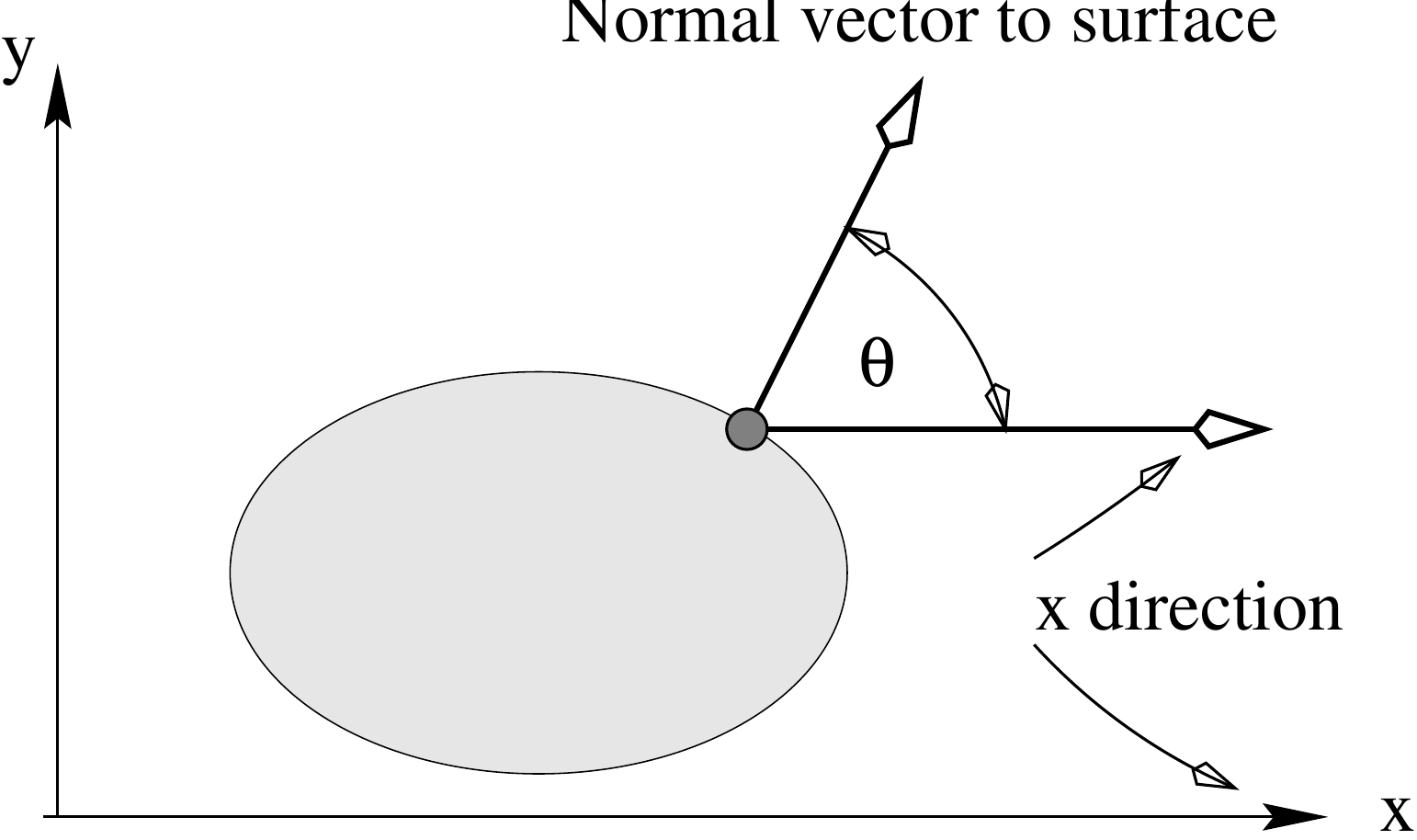}
       }
       \vspace{5.0mm}
       \centerline{
       \includegraphics[scale=0.25,angle=-90]{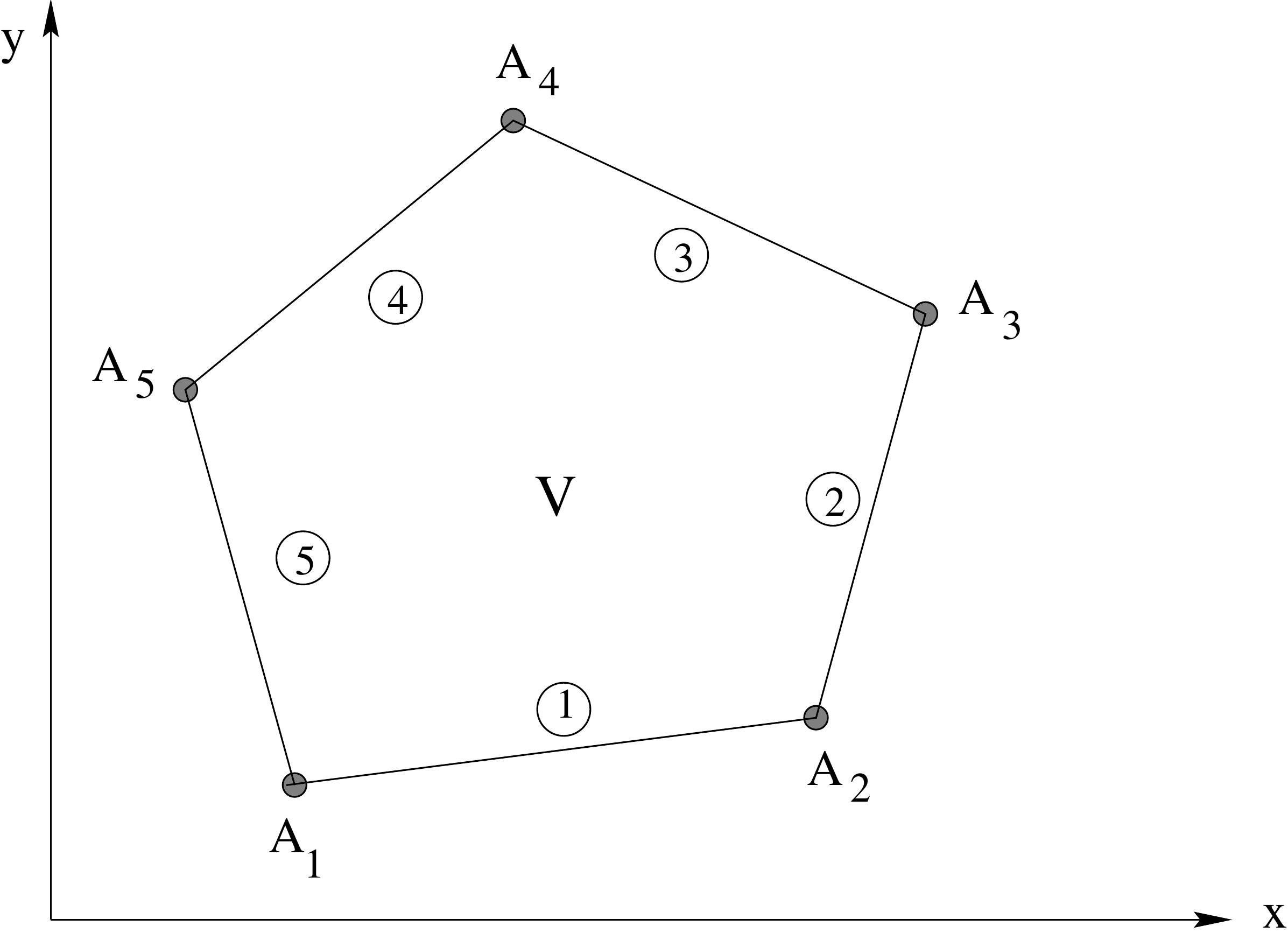}  
       }
      \caption{Control volumes, rotational invariance and finite volumes. Top frame: arbitrary control volume $V$ in Cartesian plane; the  $x$-direction is the 
      reference direction, $\theta$ is angle between the outward unit normal vector ${\bf n}$ and the 
      reference $x$-direction. Bottom frame: typical finite  
      volume in the Cartesian plane.}
\label{RotInv}
\end{figure}

As previously stated, the methods proposed in this paper combine 1D equations (\ref{swe5}) in a generic direction $s$ and the 2D equations (\ref{swe1}), locally, at a single 2D element with an arbitrary number of edges. Thus, it is convenient to recall the rotational invariance property of equations (\ref{swe1}). First, let us define an arbitrary 2D spatial control volume $V$ with boundary $\Omega$, as depicted in Fig. \ref{RotInv}, top frame. Equations (\ref{swe1}), expressed in integral form, read
\begin{equation}                              \label{rotinv1}
      \frac{\partial}{\partial t} \int_{V} {\bf Q}\, dV+ 
      \int_{\Omega}\left[\cos\theta {\bf F}({\bf Q})+\sin\theta 
      {\bf G}({\bf Q})\right]\,d \Omega  = {\bf 0}  \;.
\end{equation}
Moreover, we define the outward unit vector normal to $\Omega$, $\bf n$, as 
\begin{equation}                              \label{rotinv2}
      {\bf n}\equiv \left[n_{1},n_{2}\right] \equiv \left[\cos\theta,\sin\theta\right]\;.
\end{equation}
The top frame of Fig. \ref{RotInv} depicts the control volume $V$ in the Cartesian plane, where $x$ denotes the chosen reference direction, while the bottom frame depicts a typical computational finite volume with five vertices and five edges.
Equations (\ref{swe1}) satisfy the rotational invariance property \cite{Toro2001}
\begin{equation}                              \label{rotinv3}
      {\bf H} \equiv {\bf n} \cdot \left[{\bf F}({\bf Q}), {\bf G}({\bf Q}\right] =  
     \cos\theta {\bf F}({\bf Q}) + \sin\theta  {\bf G}({\bf Q}) = 
      {\bf T}^{-1}{\bf F}\left({\bf T}({\bf Q})\right) \; 
\end{equation}
for all vectors ${\bf Q}$ and for all real angles $\theta$, or equivalently, normal directions of $\Omega$. Here ${\bf T} = {\bf T}(\theta)$ is a rotation matrix and ${\bf T}^{-1}(\theta)$ is its inverse, given respectively as
\begin{equation}                                \label{rotinv4}
      {\bf T}=\left[\begin{array}[c]{ccc}
	     1 & 0 & 0  \\
	     0 & \cos\theta & \sin\theta \\
	     0 & -\sin\theta & \cos\theta \\              
	     \end{array}\right] \;,\hspace{2mm}
      {\bf T}^{-1}=\left[\begin{array}[c]{ccc}
	     1 & 0 & 0  \\
	     0 & \cos\theta & -\sin\theta \\
	     0 & \sin\theta & \cos\theta  \\              
	     \end{array}\right]   \;.
\end{equation}
Now, we choose a computational control volume $V_{k}$ in two-dimensional space, as shown for example in the bottom frame of Fig. \ref{RotInv}. Moreover, we define a space-time control volume $I_k = [t^{n+1},t^n] \times V_k$, over which we integrate (\ref{swe1}), yielding
\begin{equation}                                 \label{rotinv5}    
     {\bf Q}_{k}^{n+1} =  {\bf Q}_{k}^{n} -\frac{\Delta t}{ |V_{k}| }\sum_{e=1}^{N}{\cal F}_{e} \;,     		
\end{equation}
with cell averages defined as
\begin{equation}                                \label{rotinv6}    
      {\bf Q}^{n}_{k} = \frac{1}{ |V_{k}|} \int_{V_{k}} {\bf Q}(x,y,t ^{n})\, d V  \;,   		
\end{equation}
where $|V_{k}|$ denotes the volume of $V_{k}$ (or area in the 2D case), while the intercell flux for edge $e$ is
\begin{equation}                                 \label{rotinv7}    
      {\cal F}_{e}  =
      \int_{A_{e}}^{A_{e+1}} {\bf T_{e}}^{-1}{\bf F}\left({\bf T_{e}}({\bf Q})\right)\, d A \approx 
      {\cal{L}}_{e} {\bf T}_{e}^{-1} \hat{{\bf F}}_{e}   \;.    		
\end{equation}
Here ${\cal{L}}_{e}$ is the length of edge $e$, the segment $A_{e}A_{e+1}$. $\hat{{\bf F}}_{e} \approx {\bf F}\left({\bf T_{e}}({\bf Q})\right)$ is an approximation to the flux ${\bf F}$ on the edge $e$ evaluated at the rotated state ${\bf T_{e}}({\bf Q})$, where the rotation is performed in the normal direction to side $e$ through the transformation matrix ${\bf T_{e}}$. The final expression of the finite volume scheme becomes
\begin{equation}                                 \label{rotinv8}    
     {\bf Q}_{k}^{n+1} = {\bf Q}_{k}^{n} -\frac{\Delta t}{ |V_{k}| }\sum_{e=1}^{N}  {\cal{L}}_{e} {\bf T}_{e}^{-1} \hat{{\bf F}}_{e} \;.     		
\end{equation}
For completness we now illustrate the computation of the  flux $\hat{{\bf F}}_{e}$ for an arbitrary edge $e$, see  Fig.3. Conventionally, the left side $L$ of edge $e$ is always in the interior of the control volume of interest and the right side $R$ is outside. The computation of the numerical flux $\hat{{\bf F}}_{e}$, as for a first-order Godunov type method for example, involves the augmented, local one-dimensional Riemann problem in the rotated frame normal to the edge, namely
\begin{equation}                                       \label{rotinv9}
\left.
\begin{array}{ll}
      \mbox{PDEs in normal direction:}   & \partial_{t} {\bf Q}_{1D} + \partial_{s} {\bf F}_{1D}({\bf Q}_{1D}) = {\bf 0}  \;, \\
      \\
      \mbox{Rotated initial conditions:} & {\bf Q}_{1D}(s,0) = \left\{ \begin{array}{lll}
      {\bf Q}_{1D,L}  = {\bf T}_{e}({\bf Q}_{2D,L}) & \mbox{if} & s  < 0 \; , \\
      \\
      {\bf Q}_{1D,R}  = {\bf T}_{e}({\bf Q}_{2D,R}) & \mbox{if} & s > 0  \; .
      \end{array} \right.
\end{array}\right\}
\end{equation}
The steps to be followed in order to solve problem ( \ref{rotinv9}) can be summarised as:
\begin{enumerate}
\item Calculate the angle $\theta_{e}$ between the outward unit normal to edge $e$ and the fixed reference direction $x$, being positive in the counter-clock direction.

\item Calculate the corresponding rotation matrix ${\bf T}_{s}$ and its inverse from (\ref{rotinv4}).

\item Rotate left and right data as in (\ref{rotinv9}).

\item Solve the 1D Riemann problem  (\ref{rotinv9}) on rotated data and compute the corresponding flux $\hat{\bf F}_{e}$.

\item Rotate back the flux as in (\ref{rotinv3}) and multiply it by edge length to get the final intercell numerical flux for edge $e$.

\end{enumerate}
Once numerical fluxes for all edges have been calculated, element $k$ can be updated through the finite volume formula (\ref{rotinv5}). This description applies to any two-dimensional finite volume method on a general mesh, assumed here to be unstructured. More details are given in \ref{chapter_theoretical}. However, for the junction method proposed in this paper, will only use the above description at a single 2D element placed right at the junction.
\begin{figure}
       \centerline{
       \includegraphics[scale=0.35,angle=0]{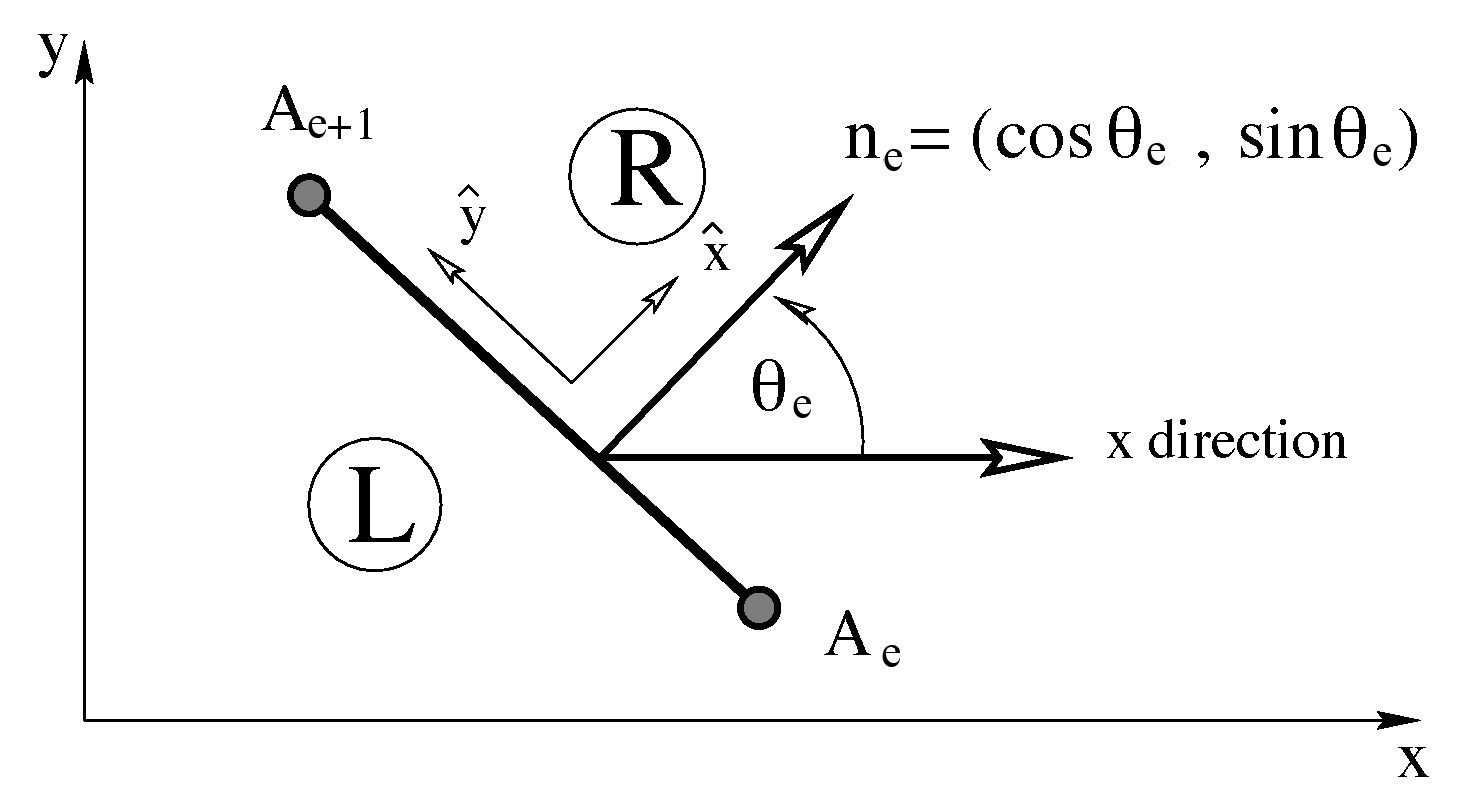}  
       }
       \vspace{3mm}
       \caption{Generic edge of a control volume $V$ in $x$-$y$ space, where by convention the left side L lies inside the control volume and the right side R lies outside. The outward normal unit vector, with respect to the fixed reference direction $x$, is depicted as well as its corresponding angle.}
\end{figure}\label{fig:genericcontrolvolumeside}

\section{A methodology for channel junctions/bifurcations} \label{sec:method}

The geometrical and numerical approaches of our proposed junction method are described below.

\subsection{The approach}

\begin{figure}[H]
\centering
\includegraphics[width=0.4\textwidth]{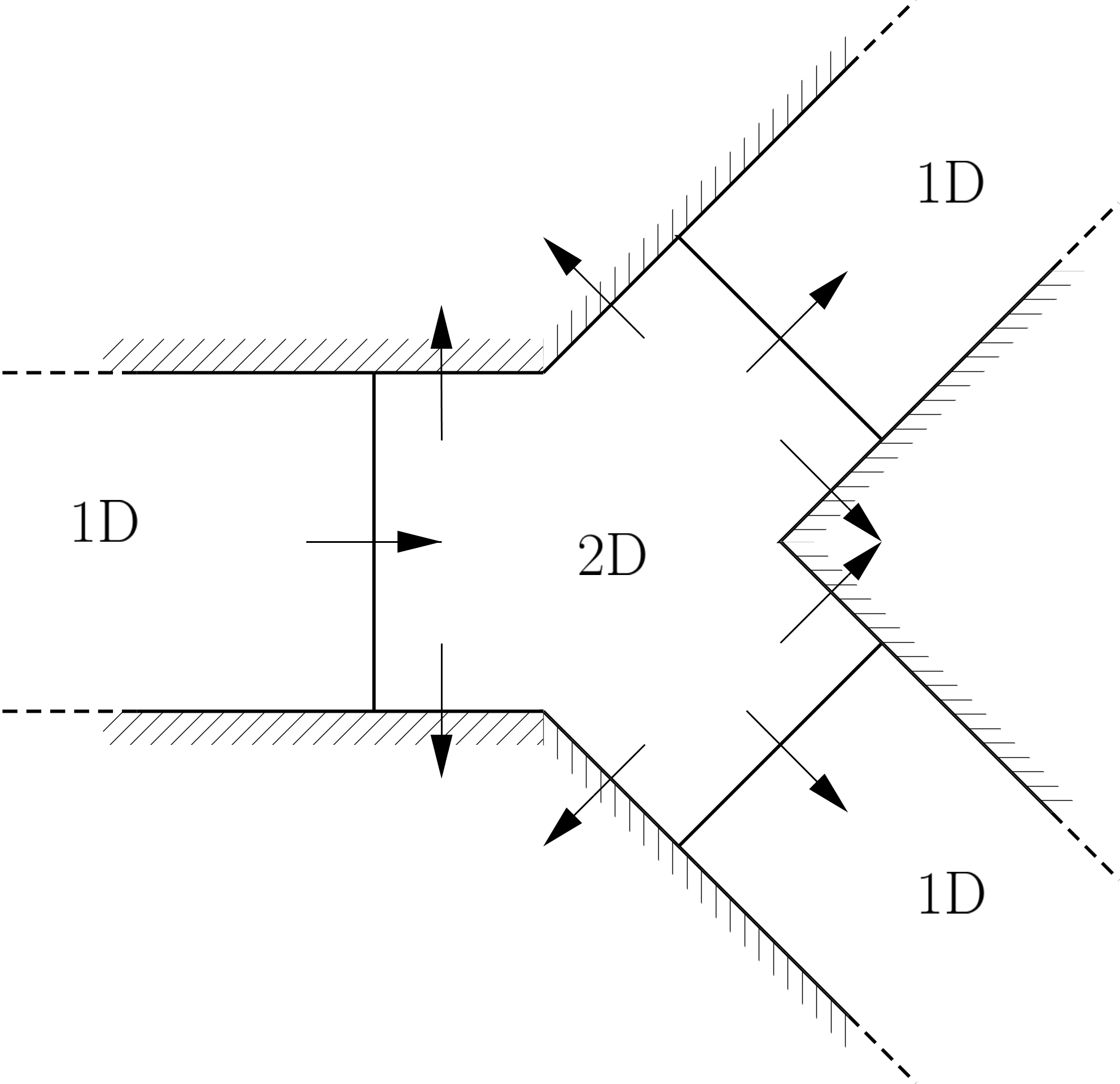}
\caption{Single 2D element at junction. The single element exchanges fluxes with all connected 1D domains. At solid walls of the rectangular cross-section channel, suitable reflective boundary conditions apply through the corersponding numerical fluxes.}
\label{methodA}
\end{figure}

In short, our method for a configuration as shown in Fig. \ref{introduzione_network} uses 1D formulations for every straight channel and a single 2D element at each junction, as depicted in Fig.  \ref{methodA}. The 2D subdomain is then linked to 1D channels through appropriate matching conditions, to be described. As previously noted, similar methods have been investigated in the past, both in hydrodynamics \cite{Miglio2005_1} and in haemodynamics \cite{Formaggia1999, FormaggiaQuarteroni2003}). Miglio et al. \cite{Miglio2005_1} used a finite element scheme and investigated only the case of  subcritical flows. In the present work we are interested in general configurations and, principally, in all possible flow regimes: subcritical, transcritical and supercritical. Our approach is independent of the particular numerical method chosen for solving the shallow water equations, but here we implement first, second and third order accurate Godunov-type finite volume methods in the ADER framework \cite{Toro:2009a}.

We remark that the choice of the shape for the 2D junction element is important and there many possible choices for fitting a single finite volume method at the junction. After investigating several possibilities we concluded that the best choice is that of a {\it junction-shaped} 2D element,  as displayed in Fig. \ref{methodA}. This 2D element protrudes into the 1D converging channels by 0.1 times the channel width, incorporating in this manner, geometrical  information on the direction of the 1D domains. Other choices for the shape of the 2D element were explored in \cite{TesiBellamoli}. The resulting method is called {\bf Method A} throughout this paper. A simple variation, called {\bf Method B}, results from the insertion of a local 2D unstructured grid composed of more than one element to represent the junction and its vicinity.

Regarding the numerical methodology for the 1D and the 2D shallow water equations we use 
Godunov-type methods with the approximate Riemann solver HLLC  \cite{HLLC}. First, second and third order accurate methods are implemented. The high-order methods follow the ADER approach \cite{ADER} with the Harten-type method to solve the generalised Riemann problem \cite{Harten}.  For background  on the ADER approach see chapters 19 and 20 of \cite{Toro:2009a}.

The time step $\Delta t$ is computed imposing a CFL condition on both the 1D and the 2D junction elements in the usual manner. Then the size $\Delta t$ is taken as the minimum among all time steps and applied to the full domain. Note that in the 2D case the maximum CFL number for stability is $CFL_{2D} = CFL_{1D}/2$. In what follows we address in more detail each of the issues arising from the coupling of 1D domains and 2D elements. 

\subsection{Computing two-dimensional and one-dimensional fluxes}

\begin{figure}[H]
\centering
\includegraphics[width=0.4\textwidth]{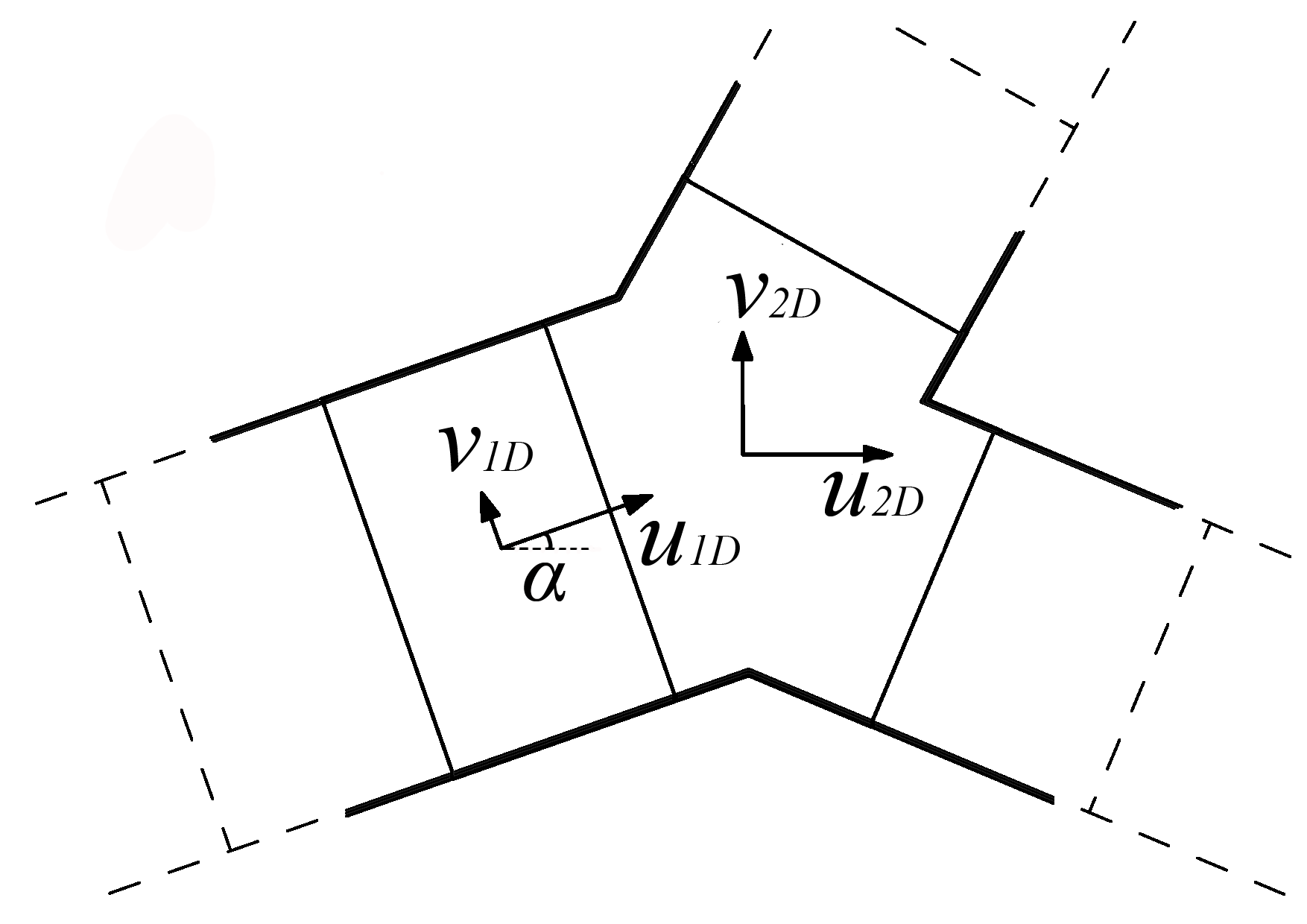}
\caption{Reference frame for the 2D element and the 1D domain of the left channel.}
\label{2DFluxesHLLC_1_MOD}
\end{figure}

We calculate the 2D fluxes  by solving the rotated 1D Riemann problem (\ref{rotinv9}) in local coordinates.  To this end we use the HLLC approximate Riemann solver \cite{HLLC}. As initial data we have $x$-velocity and $y$-velocity components in the 2D domain, where we use a global, predefined, reference system. We also have  axial velocity and transverse velocity (that is zero) in the 1D domain, where we use a local reference system, see figure \ref{2DFluxesHLLC_1_MOD}. However, for the computation of 2D fluxes we need $x$-velocity and $y$-velocity components both to the right and to the left of each edge of the element. Therefore, we need to rotate the vectors of conserved variables as follows:
\begin{equation}
{\bf Q}_{R}=\left[\begin{array}{c}
h_{2D}\\ h_{2D}u_{2D}\\ h_{2D}v_{2D}
\end{array}\right] \qquad
{\bf Q}_{L}=\left[\begin{array}{ccc}
1 & 0 & 0 \\
0 & \cos(\alpha) & -\sin(\alpha) \\
0 & \sin(\alpha) & \cos(\alpha) \\
\end{array}\right]
\left[\begin{array}{c}
h_{1D}\\ h_{1D}u_{1D}\\ h_{1D}v_{1D}
\end{array}\right]\;,
\end{equation}
where $u_{1D}$, $v_{1D}$ and $h_{1D}$, are the variables in the 1D domain, while $u_{2D}$, $v_{2D}$ and $h_{2D}$ denote the variables in the 2D domain, as depicted in Fig. \ref{2DFluxesHLLC_1_MOD}. For a second or higher order scheme, these data values result from reconstructed polynomials evaluated at the edges. Once the variables are available, we can apply the classical HLLC solver \cite{HLLC} as described in \ref{chapter_theoretical}. For each channel and for each side of the junction there is a reference angle which we call $\alpha$. 
Reflective boundary conditions are set on the remaining edges of the junction-shaped element, giving rise to symmetric Riemann problems, see \cite{Toro:2009a} for details. 
With regard to the 1D channel on the left side of Fig. \ref{2DFluxesHLLC_1_MOD}, the problem is inverted; we need  axial velocity and transversal velocity both to the rigth and to the left of each edge of the 1D cell. Vectors of conserved variables become:
\begin{equation}
{\bf Q}_{L}=\left[\begin{array}{c}
h_{1D}\\ h_{1D}u_{1D}\\ h_{1D}v_{1D}
\end{array}\right] \qquad
{\bf Q}_{R}=\left[\begin{array}{ccc}
1 & 0 & 0 \\
0 & \cos(\alpha) & \sin(\alpha) \\
0 & -\sin(\alpha) & \cos(\alpha) \\
\end{array}\right]
\left[\begin{array}{c}
h_{2D}\\ h_{2D}u_{2D}\\ h_{2D}v_{2D}
\end{array}\right]\;.
\end{equation}

\subsection{Dealing with transverse velocity in 1D domains}

Obviously, in all 1D elements we assume 1D motion and thus the transverse velocity component is zero.  
However, a problem arises at the element of a 1D channel adjacent to a 2D junction element, since we might end up with a non-zero transverse velocity. In  the 2D elements, at time $t^{n}$ we generally will have two non-zero velocity components, and consequently we could obtain a non-zero transversal velocity also in the 1D element at time $t^{n+1}$ due to the 2D flux at the edge of the 2D element. To deal with this difficulty we have considered two approaches. One possibility is to simply set to zero the transverse velocity and take as 1D axial velocity the normal velocity component. The second option, which we prefer,  is to calculate the axial 1D velocity as
\begin{equation}
  u_{1D}^{n+1}=\text{sign}(u_{1D}^{n+1})\sqrt{(u_{2D}^{n+1})^2+(v_{2D}^{n+1})^2}\;.
\end{equation}
This means that the 2D velocity vector has been rotated in the direction of the 1D channel and consequently the tranverse velocity component is zero. Inevitably, in both approaches, momentum balance at the 1D elements adjacent to the 2D element is effectively altered, even though in the 2D element momentum balance is strictly satisfied.

\subsection{Spatial reconstruction for high-order accuracy}

\begin{figure}[H]
\centering
\subfigure[Cells used for 2D reconstruction.]{\label{2DFluxesHARTEN_1_MOD}
\includegraphics[width=0.4\textwidth]{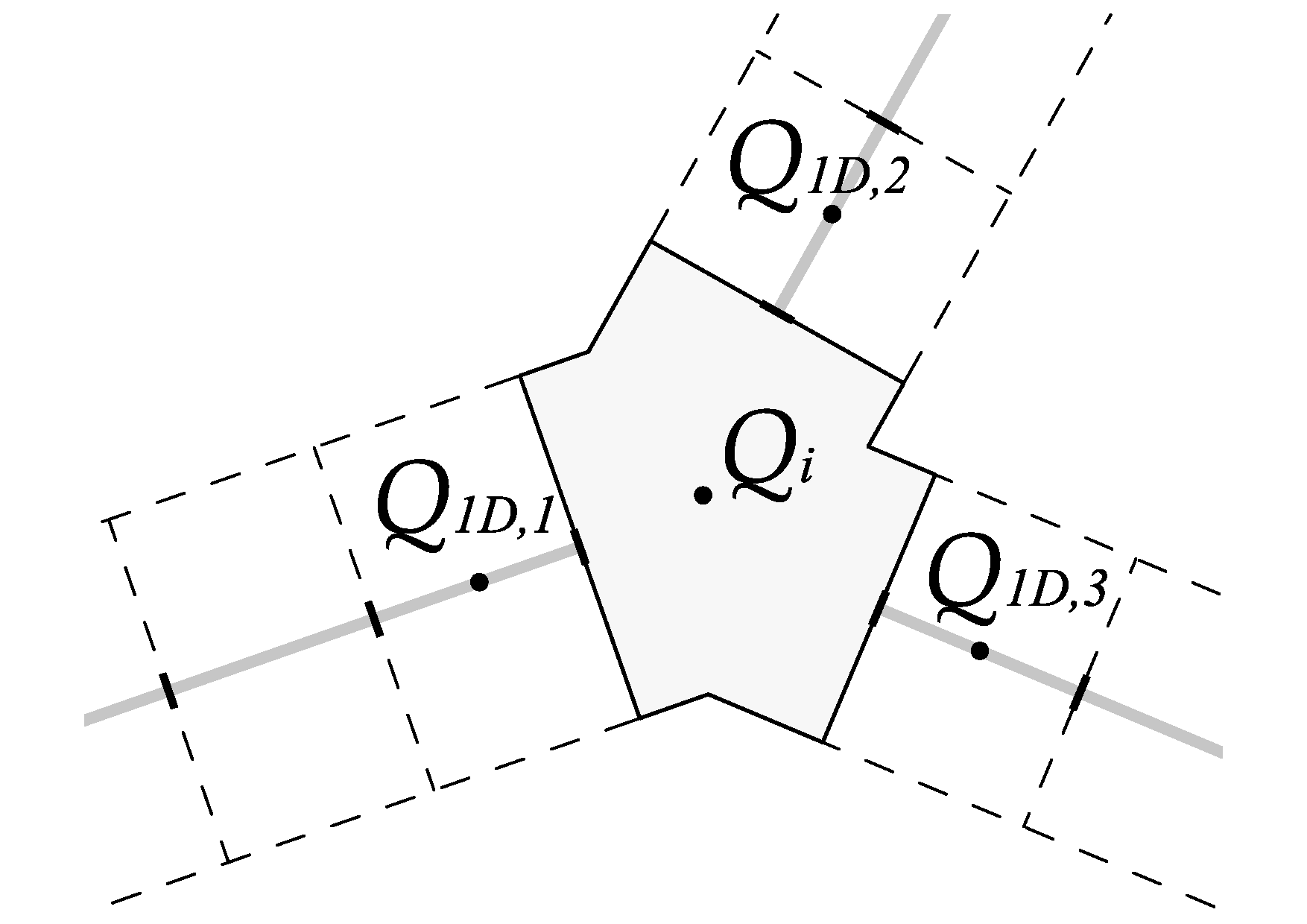}}\qquad
\subfigure[One-dimensional reconstruction.]{\label{1DFluxesHARTEN_2_MOD}
\includegraphics[width=0.4\textwidth]{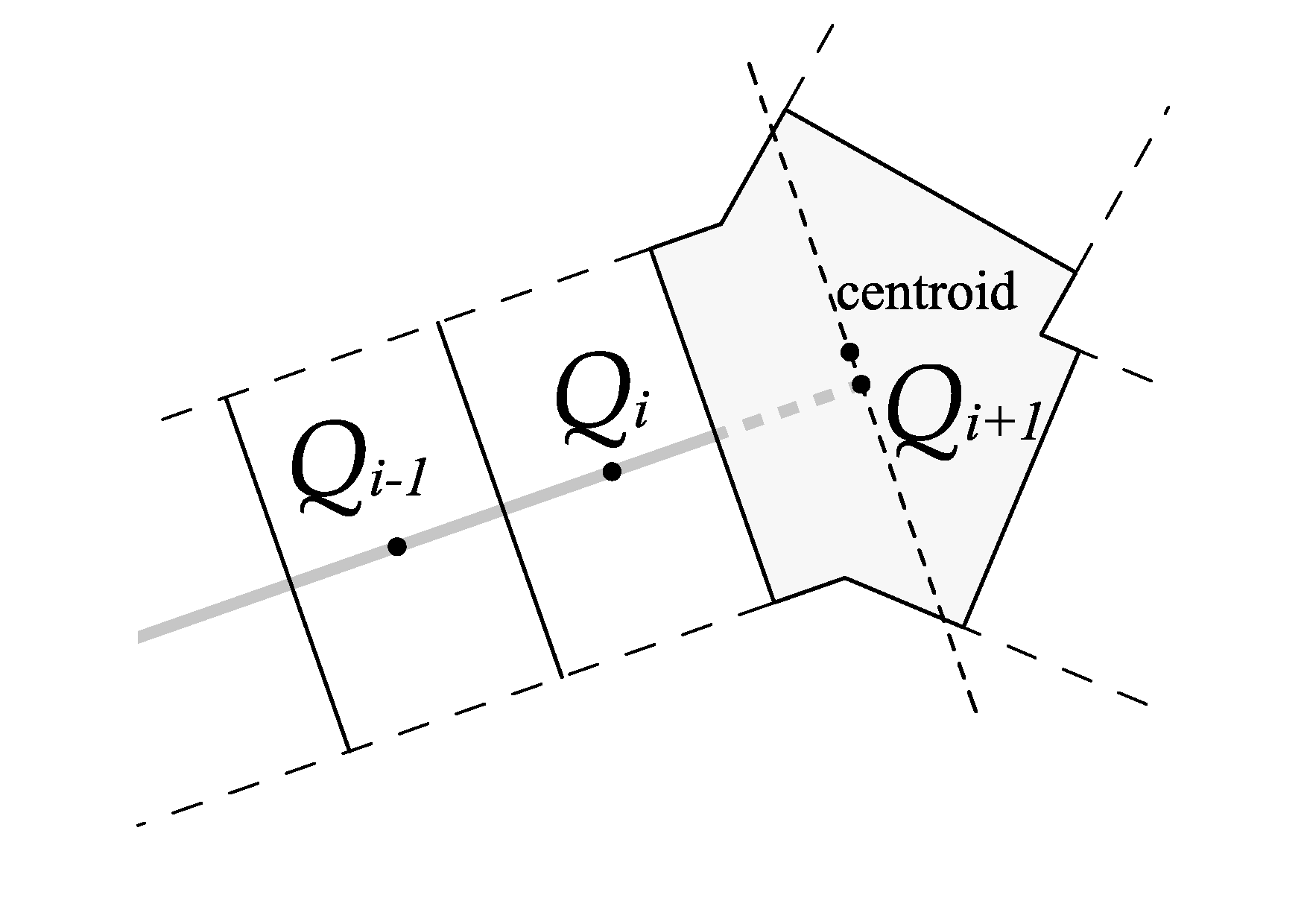}}
\caption{Illustration and notation for the spatial reconstruction in 2D (a) and 1D (b).}
\end{figure}

For 1D cells adjacent to 2D cells we perform a modified version of the spatial reconstruction described in \ref{chapter_theoretical}, by projecting the distance between the centroid of the 2D element and the centre of the 1D cell along the normal to the boundary. See figure \ref{1DFluxesHARTEN_2_MOD}.
Concerning 2D elements, particular attention must be paid to the reconstruction process. As in the 1D case, at any given time level $n$ one has a set of constant volume averages that are approximations to integral averages within each finite volume. For a second-order scheme, we need to approximate the solution in the 2D element with a first order polynomial. To this end we need three equations, for which we consider the three neighbouring 1D cells, as shown in Fig. \ref{2DFluxesHARTEN_1_MOD}. We do not use fictious elements near reflective boundaries for the reconstruction.

The 1D reconstruction delivers a slope in the axial direction, while 2D reconstructions results in slopes in $x$- and $y$-directions. When passing from 1D to 2D or vice versa we need to transform the first into the second, so we have to rotate not only the vector of conserved variables, but also the gradients. In fact, in the 1D domain we have $\partial_n u$ and $\partial_n v$, but to apply Harten's approach to solve the generalized Riemann problem we need $\partial_x U$, $\partial_y U$, $\partial_x V$ and $\partial_y V$ (being $u$ and $v$ velocities in axial and transversal directions, and $U$ and $V$ velocity components in $x$ and $y$ directions). These slopes can be calculated as
\begin{equation}
\left(\begin{array}{c}
\partial_x u\\
\partial_y u
\end{array}\right)=
\left(\begin{array}{c}
\cos\alpha\\
\sin\alpha
\end{array}\right)\partial_n u\;, \qquad\quad
\left(\begin{array}{c}
\partial_x v\\
\partial_y v
\end{array}\right)=
\left(\begin{array}{c}
\cos\alpha\\
\sin\alpha
\end{array}\right)\partial_n v \;
\end{equation}
and
\begin{equation}
\begin{array}{c}
\left(\begin{array}{c}
\partial_x U\\
\partial_x V
\end{array}\right)=
\left[\begin{array}{cc}
\cos\alpha & -\sin\alpha\\
\sin\alpha & \cos\alpha
\end{array}\right]
\left(\begin{array}{c}
\partial_x u\\
\partial_x v
\end{array}\right) \;,         \\
\\
\left(\begin{array}{c}
\partial_y U\\
\partial_y V
\end{array}\right)=
\left[\begin{array}{cc}
\cos\alpha & -\sin\alpha\\
\sin\alpha & \cos\alpha
\end{array}\right]
\left(\begin{array}{c}
\partial_y u\\
\partial_y v
\end{array}\right) \;.
\end{array}
\end{equation}
In the next section we assess the performance of the junction methods presented in this paper using a comprehensive suite of test problems, comparing results to 2D reference solutions obtained from an unstructured 2D second-order method discribed in \ref{chapter_theoretical}.

\section{Test problems and assesment of the methods} \label{sec:results}

\begin{figure}[H]
\centering
\includegraphics[width=0.3\textwidth]{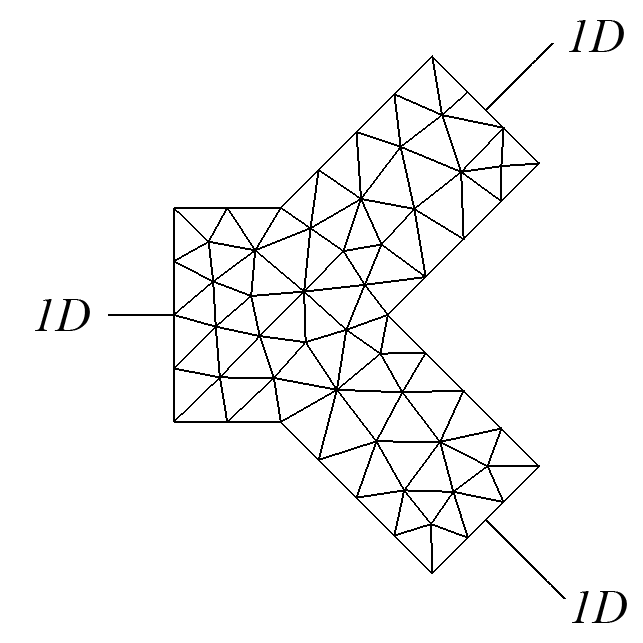}
\caption{Example of a local 2D grid used in the vicinity of the junction region in {\bf Method B}.}
\label{methodB}
\end{figure}

In this paper we consider three methods to deal with junctions in the context of shallow-water channels, {\bf Method A} being our main contribution, in which a single 2D element is inserted in each junction. {\bf Method B} generalises {\bf method A} by inserting a local 2D unstructured grid in the vicinity of each junction, see Fig. \ref{methodB}. The third method considered for comparison is the method proposed by Peir\'o, Sherwin, Formaggia and Parker \cite{SherwinFormaggia,Sherwin2003}, which in this paper will be called {\bf PSFP method}. This method is summarised in \ref{section_1Dexisting}. Solutions from all 3 methods are compared to reference 2D solutions. All methods have been implemented to second-order accuracy in both space and time. Here we present results for five tests. For additional tests see \cite{TesiBellamoli}.

\subsection{Single-junction test problems}


\noindent{\bf Test 1: Subcritical wave in channel with a $90^\circ$ bifurcation.}\label{testsub90}\\

\begin{figure}[H]
\centering
\addtocounter{subfigure}{-1}
\subfigure{\raisebox{10mm}{\includegraphics[width=0.15\textwidth]{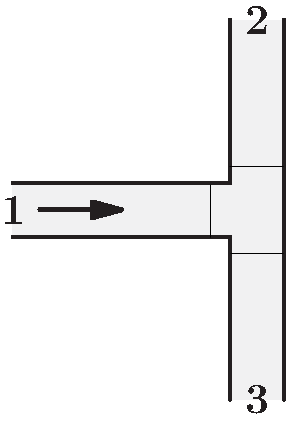}}}
\hspace{0.05\textwidth}
\subfigure[Channel 1]{\label{junct90onda_t8s_ch1_ordine}
\includegraphics[width=0.65\textwidth]{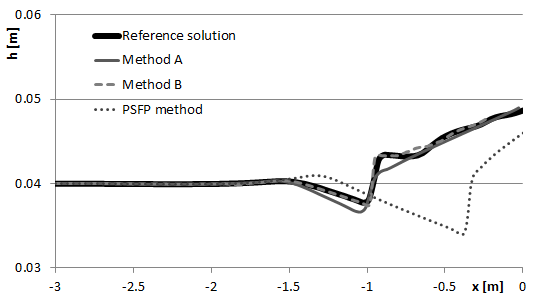}}
\addtocounter{subfigure}{-1}
\subfigure{\raisebox{10mm}{\includegraphics[width=0.16\textwidth]{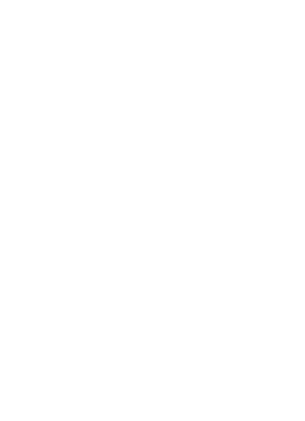}}}
\hspace{0.05\textwidth}
\subfigure[Channel 2]{\label{junct90onda_t8s_ch2_ordine}
\includegraphics[width=0.65\textwidth]{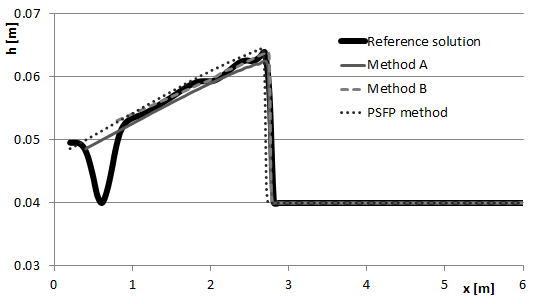}}
\caption{Test 1: Subcritical wave. Water height at time $t=8\,s$. } \label{junct90onda}
\end{figure}

In this test we consider a channel configuration as shown on the left of Fig. \ref{junct90onda}. We impose a subcritical wave ($Fr_{max}\simeq 0.4$) that gradually steepens up and becomes a shock wave just after a $90^\circ$ bifurcation. Results are shown in Fig. \ref{junct90onda}.  Method A and B give very satisfactory results, as compared to the 2D reference solution, for channel 1, while the PSFP method gives rather inaccurate results. For channel 2 all three methods give quite similar results, methods A and B being slightly more accurate than the PSFP method.

\noindent{\bf Test 2: Subcritical wave in a channel with a $90^\circ$ asymmetrical bifurcation.}\\

\begin{figure}[H]
\centering
\addtocounter{subfigure}{-1}
\subfigure{\raisebox{10mm}{\includegraphics[width=0.3\textwidth]{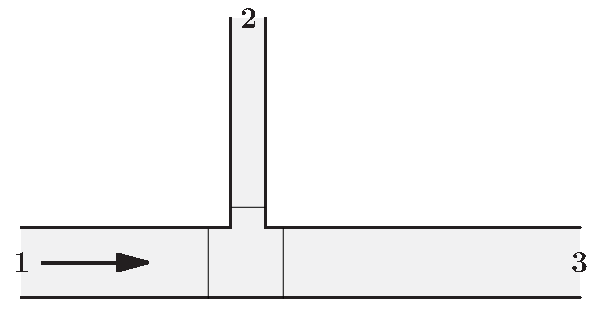}}}
\hspace{0.03\textwidth}
\subfigure[Channel 1]{\label{junct90ASonda_t8s_ch1}
\includegraphics[width=0.65\textwidth]{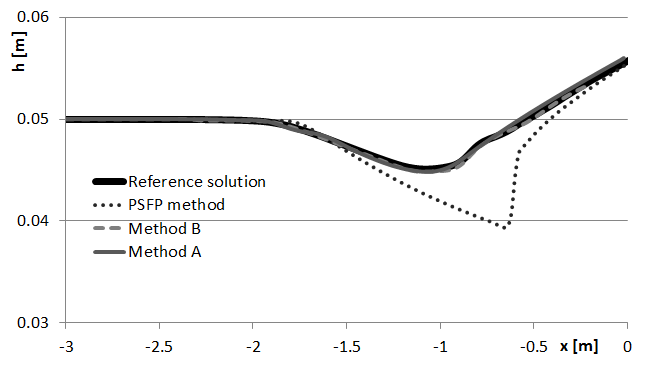}}
\addtocounter{subfigure}{-1}
\subfigure{\raisebox{10mm}{\includegraphics[width=0.3\textwidth]{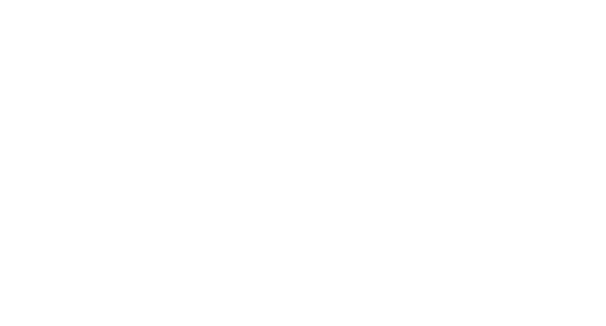}}}
\hspace{0.03\textwidth}
\subfigure[Channel 2]{\label{junct90ASonda_t8s_ch2}
\includegraphics[width=0.65\textwidth]{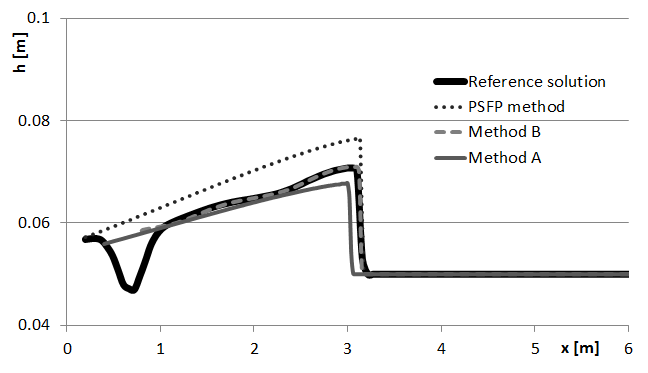}}
\caption{Test 2: Subcritical wave (asymmetrical case): Water height  at time $t=8\,s$.} \label{junct90ASonda}
\end{figure}

In this test we consider an asymmetrical  channel configuration as shown on the left of Fig. \ref{junct90ASonda}. As for the previous test, methods A and B give very satisfactory results as compared to the reference 2D solution, outperforming the PSFP method.

\noindent{\bf Test 3: Shock wave in a channel with a $45^\circ$ bifurcation.}\\

In this test we consider a  channel configuration as shown on the left of Fig. \ref{junct90shocksuper}. 
From channel 1 we send a shock with Froude number $Fr=0.75$. 
Results are shown in Fig. \ref{junct90shocksuper}. It is seen that the performance of  methods A and B is very satisfactory, as far as the shock wave is concerned. The PSFP method did not run for this test.
\begin{figure}[H]
\centering
\addtocounter{subfigure}{-1}
\subfigure{\raisebox{5mm}{\includegraphics[width=0.2\textwidth]{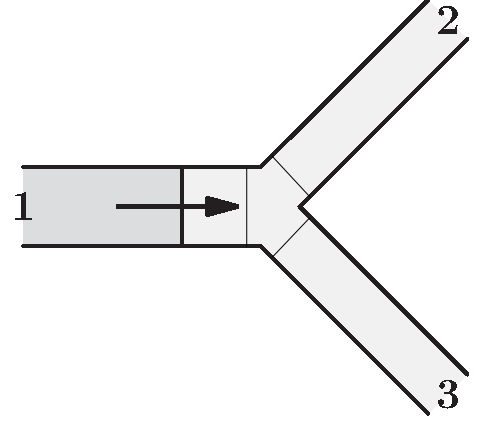}}}
\hspace{0.05\textwidth}
\subfigure[Channel 2]{\label{junct90shocksuper_t2s_ch2}
\includegraphics[width=0.65\textwidth]{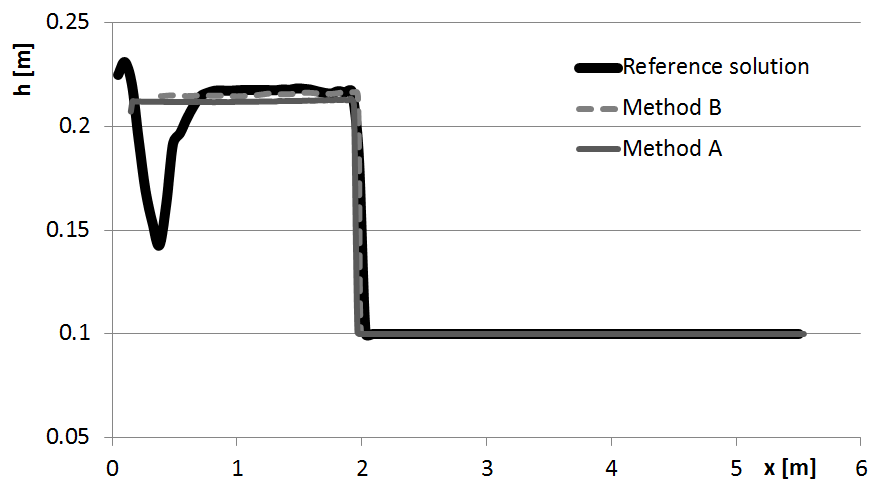}}
\caption{Test 3: Supercritical shock wave ($45^\circ$). Water height at time $t=2\,s$.} \label{junct90shocksuper}
\end{figure}

\noindent{\bf Test 4: Supercritical shock wave in a channel with a $90^\circ$ bifurcation.}\\

Finally we test our methods with a severe problem: a supercritical shock of Froude number $Fr=1.135$. Results are shown in Fig. \ref{junct90shocksuper}. Results obtained with method B are again very satisfactory, thanks to the local 2D grid. On the other hand, results obtained with method A are less accurate than that obtained in the previous case, because of the severity of the test. Again, the PSFP method did not run for this test.
\begin{figure}[H]
\centering
\addtocounter{subfigure}{-1}
\subfigure{\raisebox{5mm}{\includegraphics[width=0.15\textwidth]{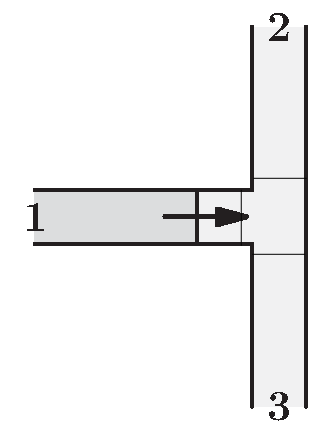}}}
\hspace{0.05\textwidth}
\subfigure[Channel 2]{\label{junct90shocksuper_t2s_ch2}
\includegraphics[width=0.65\textwidth]{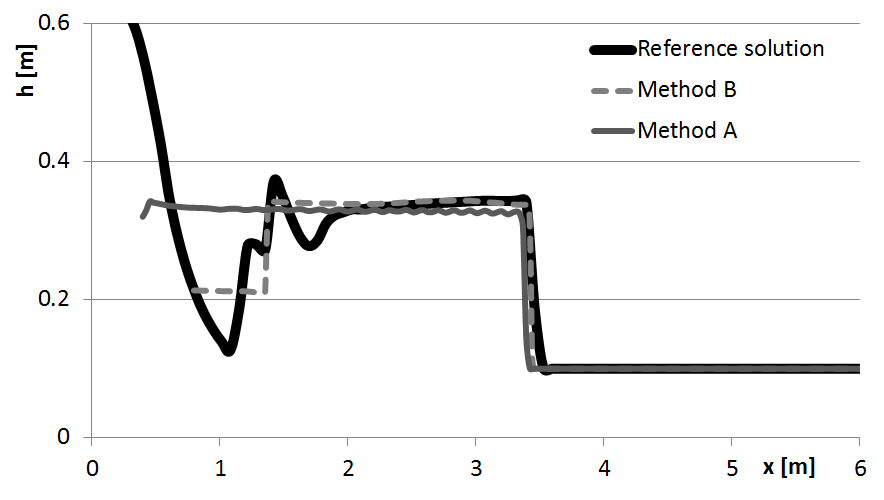}}
\caption{Test 4: Supercritical shock wave ($90^\circ$ bifurcation). Water height at time $t=2\,s$.} \label{junct90shocksuper}
\end{figure}

\subsection{\bf Test 5: the CADAM test problem.}

\begin{figure}[H]
\centering
\subfigure[Method B]{\label{cadam_grid}
\includegraphics[width=0.48\textwidth]{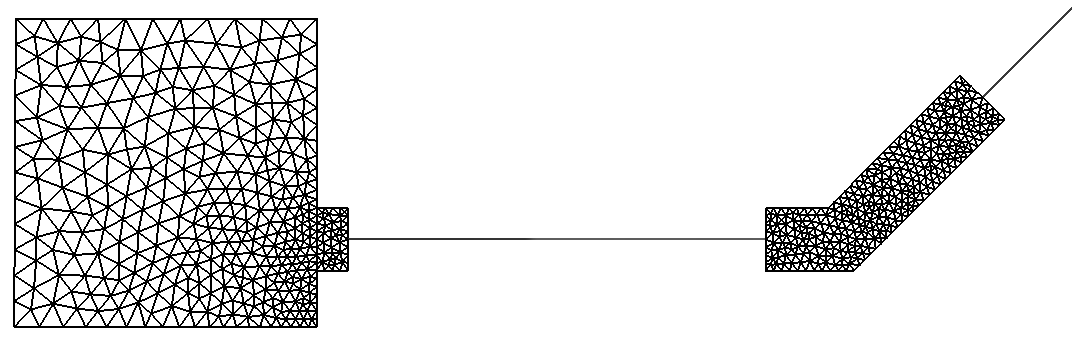}}
\subfigure[Method A]{\label{cadam_element}
\includegraphics[width=0.48\textwidth]{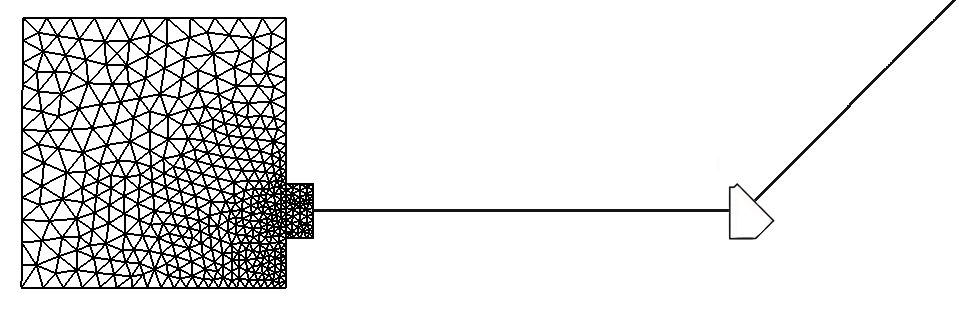}}
\caption{Test 5: the CADAM test problem. 2D and 1D domains used for numerical simulation of CADAM test 1.}
\label{cadam}
\end{figure}

\begin{figure}[H]
\centering
\subfigure[Gauge 2]{
\includegraphics[width=0.45\textwidth]{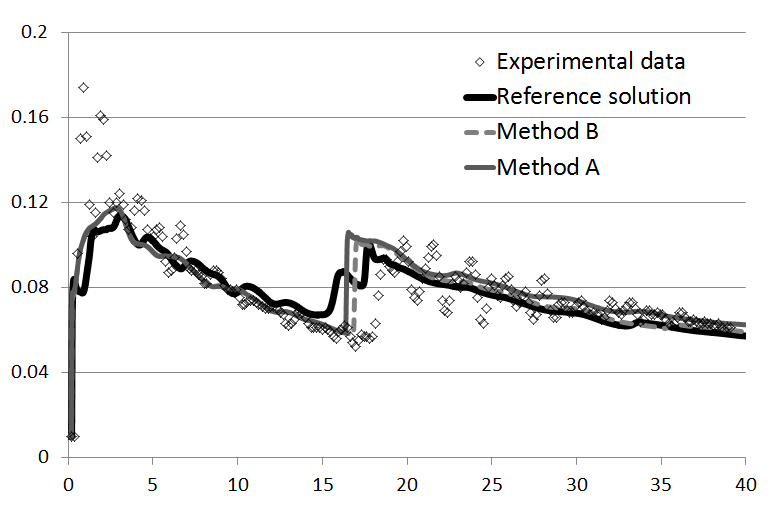}}
\subfigure[Gauge 3]{
\includegraphics[width=0.45\textwidth]{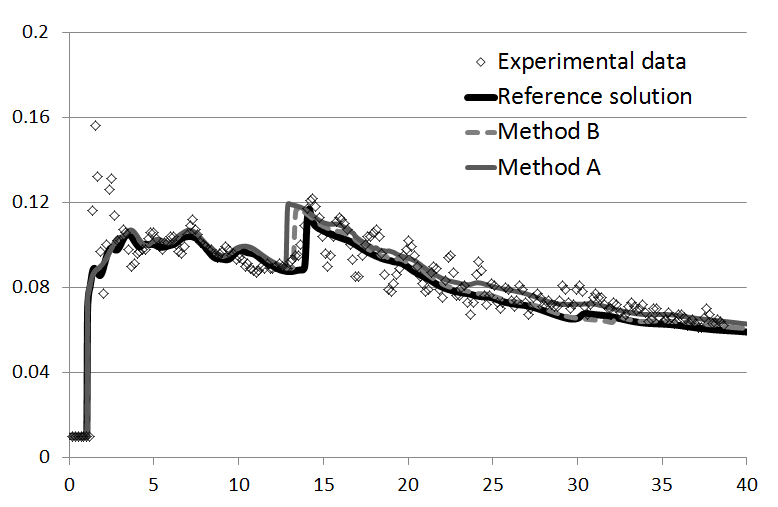}}
\subfigure[Gauge 4]{
\includegraphics[width=0.45\textwidth]{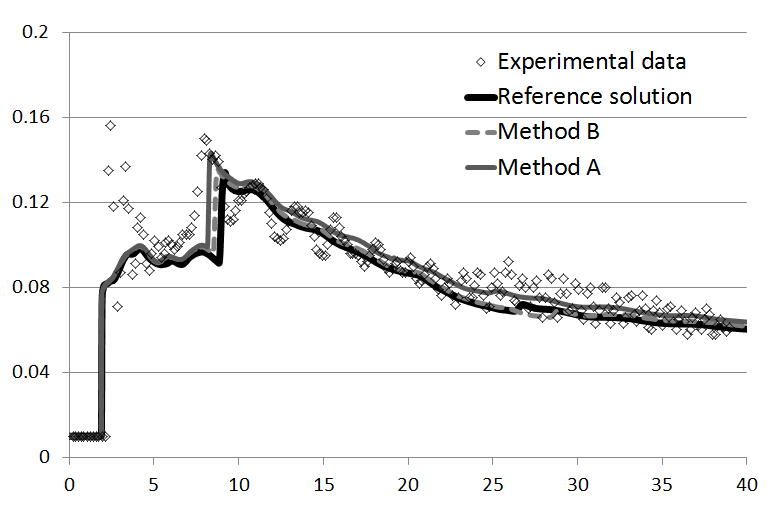}}
\subfigure[Gauge 5]{
\includegraphics[width=0.45\textwidth]{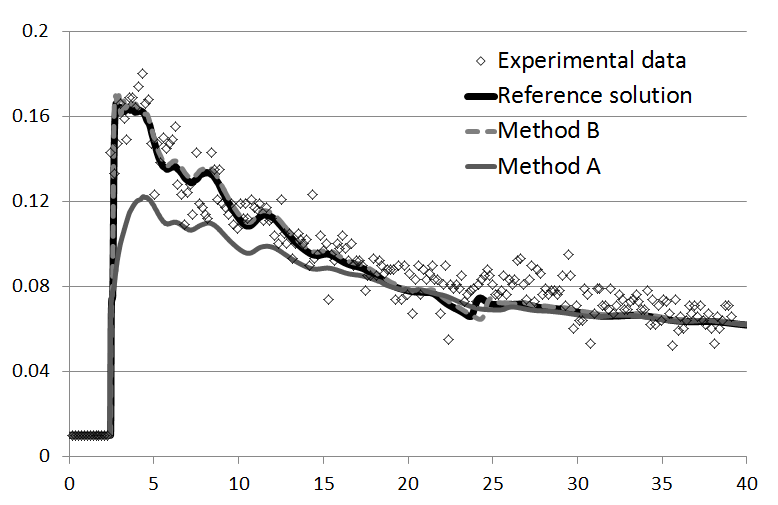}}
\subfigure[Gauge 6]{
\includegraphics[width=0.45\textwidth]{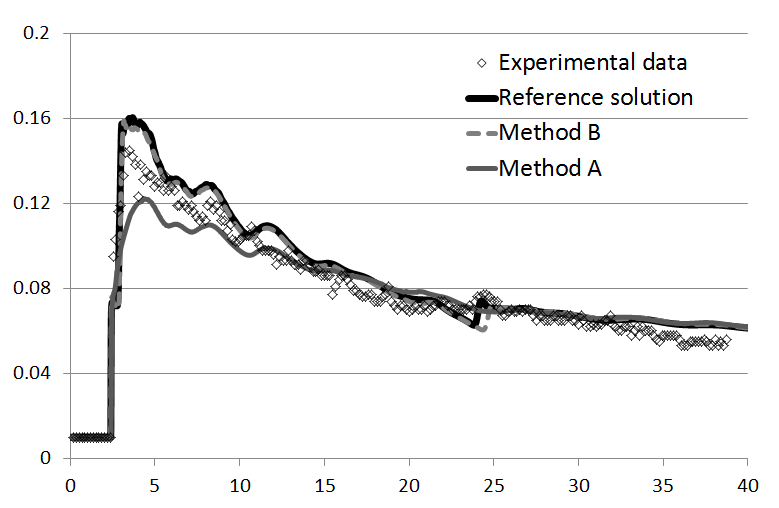}}
\subfigure[Gauge 7]{
\includegraphics[width=0.45\textwidth]{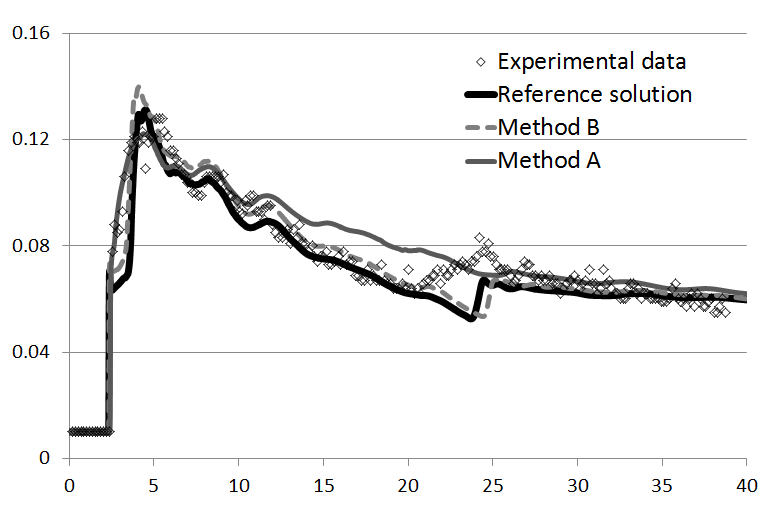}}
\subfigure[Gauge 8]{
\includegraphics[width=0.45\textwidth]{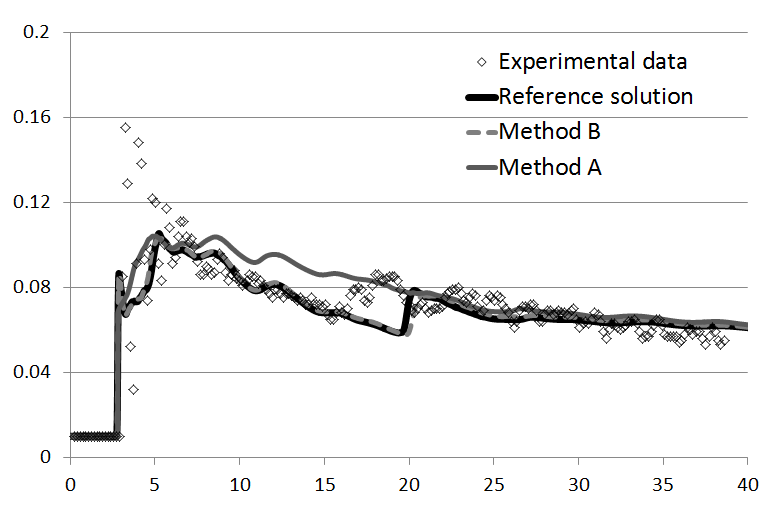}}
\subfigure[Gauge 9]{
\includegraphics[width=0.45\textwidth]{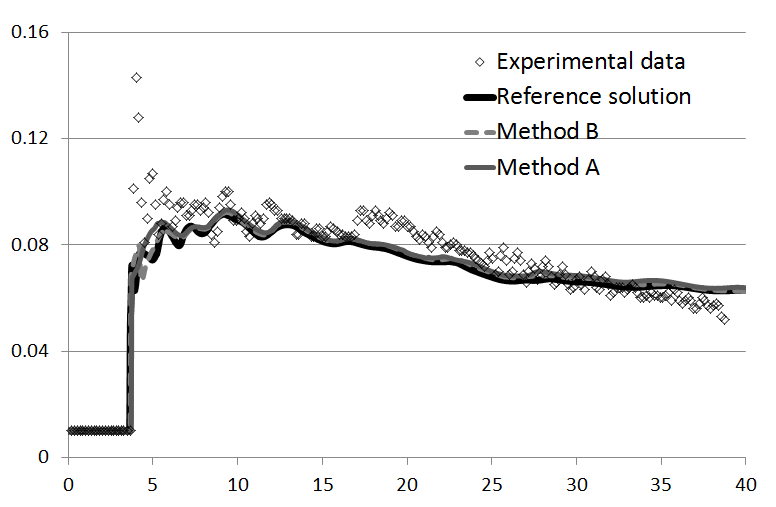}}
\caption{Test 5: the CADAM test problem. Computed free-surface elevation [meters] in time [seconds] and experimental measurements. Gauges 2 to 9 are the points of measurement used in the experimental test  \cite{Morris}.}
\label{cadam_coupled_source}
\end{figure}

In this section we apply the methods to the CADAM test 1 (CADAM, Concerted Action on Dam-Break Modelling, 1996-1999), for which experimental measurements are available as well as numerous numerical simulations. For a full description of the test see \cite{Morris}. The geometrical configuration is depicted in Fig. \ref{cadam} in which a 2D reservoir is connected to a straight channel with a $45^\circ$ bend. Figs. \ref{cadam_grid} and \ref{cadam_element} show how the the $45^\circ$ bend was treated for methods A and B. In both cases the reservoir is discretised with a 2D unstructured mesh, while for the $45^\circ$ bend method B inserts a local 2D grid in the vicinity of the bend while method A considers a single 2D element. 
In the CADAM experiment, measuring gauges 5 to 7 are placed around the bend, where the motion of the fluid is more complex. Gauges 2, 3, 4 and 9 are placed along the straight channels. For full details see \cite{Morris}.  Numerical results and experimental measurements are all displayed in Fig \ref{cadam_coupled_source}. Results obtained with the methods proposed in this paper compare satisfactorily to measurements. The flow is supercritical, so the PSFP method did not run for this test. 

\subsection{Test 6: A multiple-channel network}

In this section we assess the performance of the various methods for the case of a multiple-channel network involving 16 junctions and 25 branches; see Figs. 14-16. We considered two cases, an incident subcritical wave and an incident supercritical shock. For the sake of simplicity we set the bed slope and the friction to zero.  Solutions are computed with all three approximate junction methods considered, except for the supercritical shock case for which only methods A and B are used.  For this test, due to the complexity of the situation with many shock wave reflections and wave interaction, for method A we use a 2D coarse grid inside the four junctions in the corners (see Fig. \ref{rete_elementgrid}), where the flow is very complex due to large variations in angles and the large space occupied by the junction. Results will be shown at the eight positions shown in Fig. \ref{rete3}. 
\begin{figure}[H]
\centering
\includegraphics[width=0.85\textwidth]{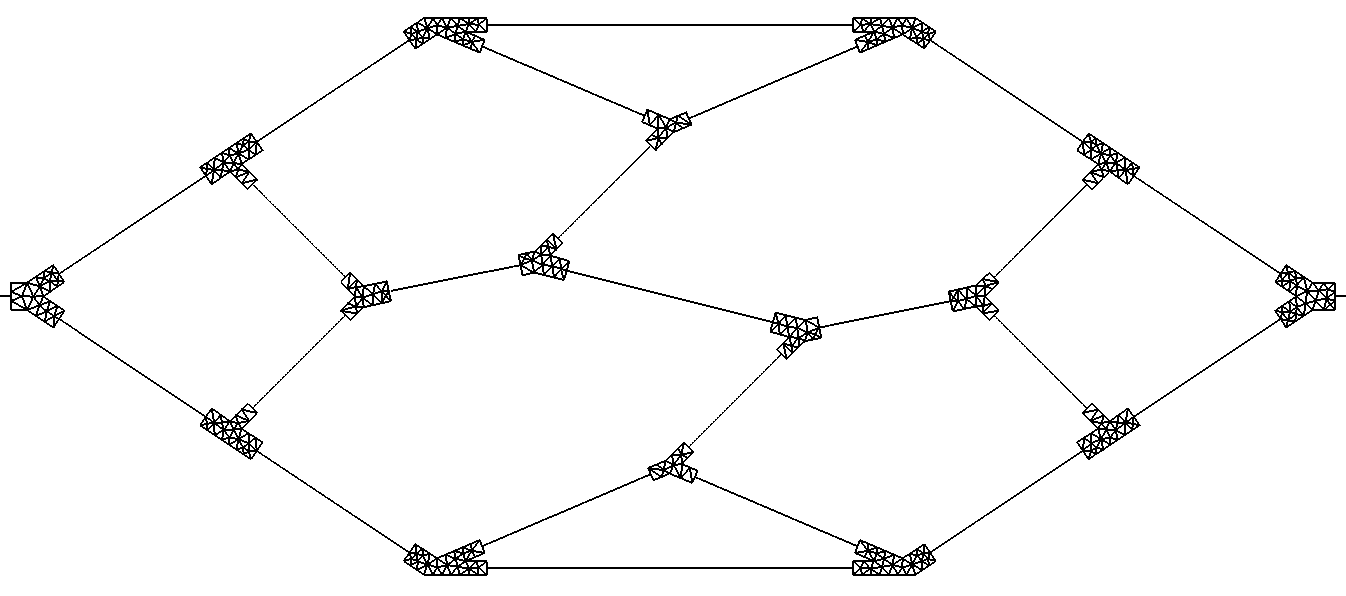}
\caption{Test 6: A multiple-channel network. Configuration for method B. Two-dimensional grids in the vicinity of junctions.}
\label{rete_grid}
\end{figure}
\begin{figure}[H]
\centering
\includegraphics[width=0.85\textwidth]{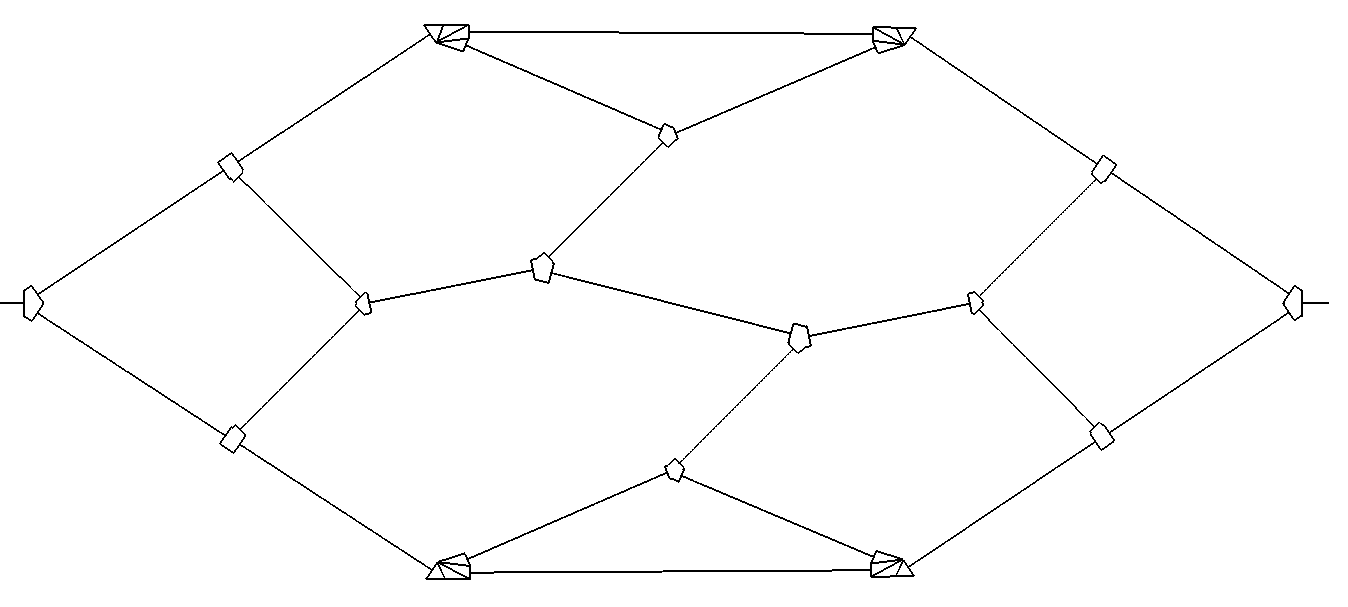}
\caption{Test 6: A multiple-channel network. Method A. 2D single elements in most junctions, modified at shown four corner junctions.}
\label{rete_elementgrid}
\end{figure}
For the subcritical wave case, computed results are displayed in Figs. \ref{rete_onda_grafici1} to \ref{rete_onda_grafici8}. All three approximate junction methods run and are compared to the reference 2D solution. Methods A and B are seen to be very accurate; all three methods give very similar results for the arrival phase of the wave but differ at later times. For the supercritical shock wave case, computed results are displayed in Figs. \ref{rete_shock_grafici1} to \ref{rete_shock_grafici8}. For this case the PSFP method did not run. Not surprisingly, it is seen that method B hardly differs from the reference 2D solution but the simpler method A is also seen to be very accurate. As expected, the larger discrepancies between method A and the reference solution are seen in wave arrival times.  Results at point 8 were expected to show the largest errors, as waves must transverse the full complex network, with multiples shock waves and complex interactions, and yet the end results at position 8 are satisfactory.

%
\begin{figure}[H]
\centering
\includegraphics[width=0.85\textwidth]{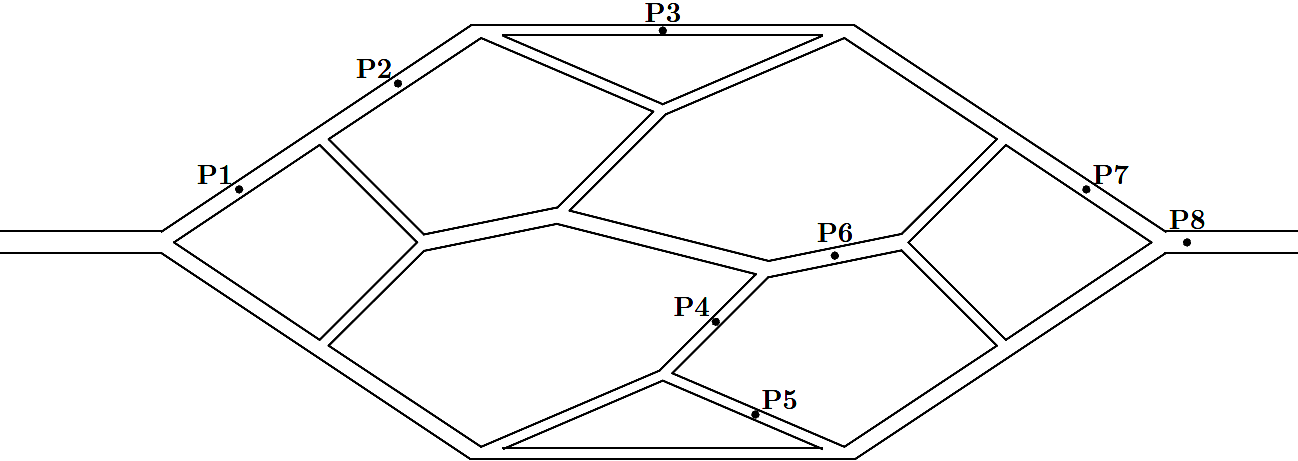}
\caption{Test 6: A multiple-channel network. Points of the network where the free surface elevation is recorded and then reported in figures \ref{rete_onda_grafici1} to \ref{rete_shock_grafici8}.}
\label{rete3}
\end{figure}
\begin{figure}[H]
\centering
\includegraphics[width=0.7\textwidth]{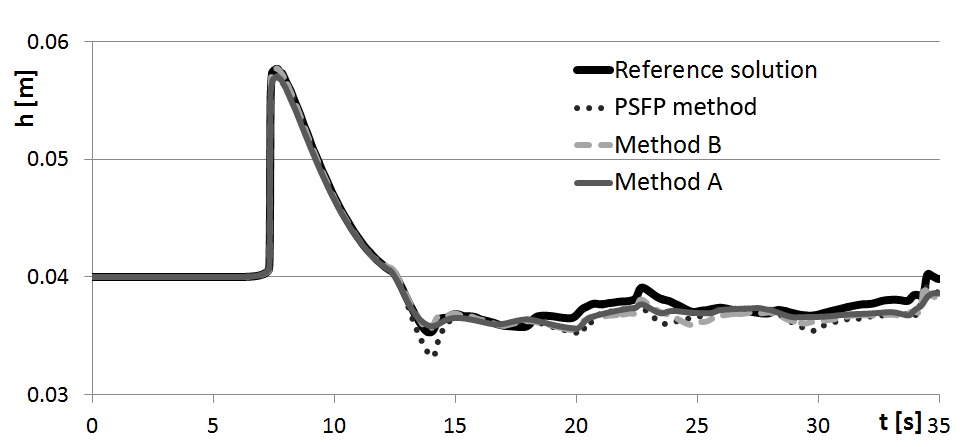}
\caption{Test 6 (subcritical wave): A multiple-channel network. Computed free-surface elevation [m] in time [s] for Point 1.}
\label{rete_onda_grafici1}
\end{figure}
\begin{figure}[H]
\centering
\includegraphics[width=0.7\textwidth]{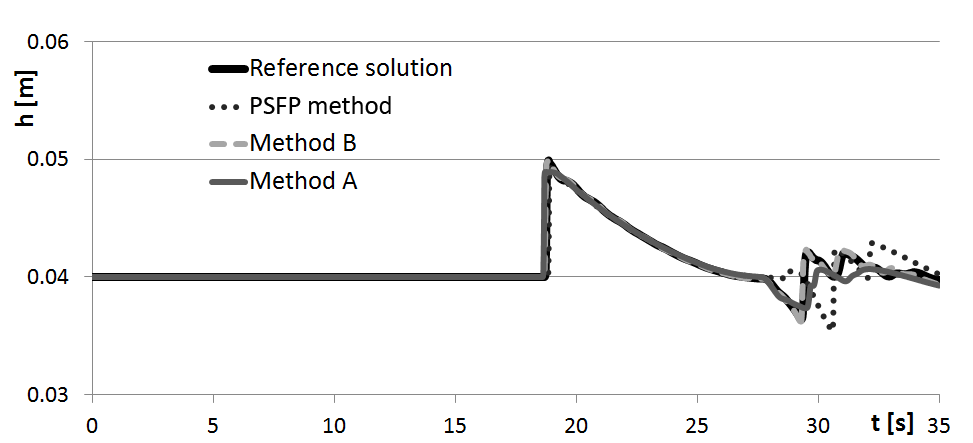}
\caption{Test 6 (subcritical wave): Computed free-surface elevation [m] in time [s] for Point 3.}
\label{rete_onda_grafici3}
\end{figure}
\begin{figure}[H]
\centering
\includegraphics[width=0.7\textwidth]{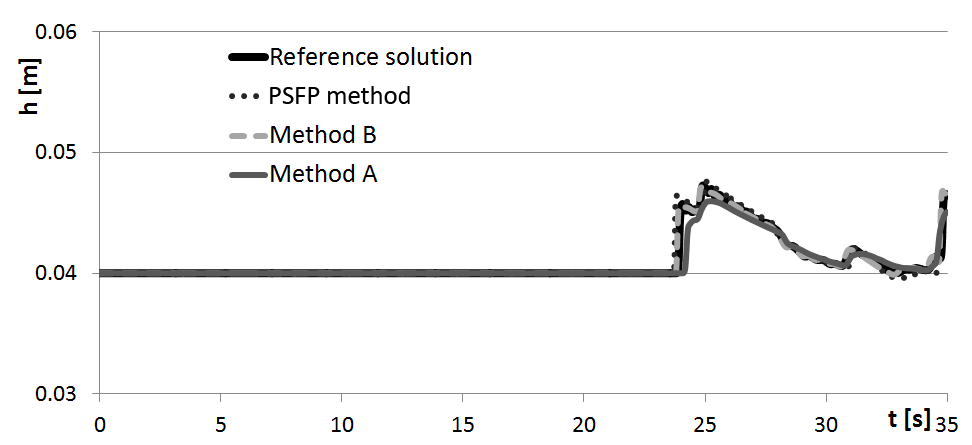}
\caption{Test 6 (subcritical wave): Computed free-surface elevation [m] in time [s] for Point 6.}
\label{rete_onda_grafici6}
\end{figure}
\begin{figure}[H]
\centering
\includegraphics[width=0.7\textwidth]{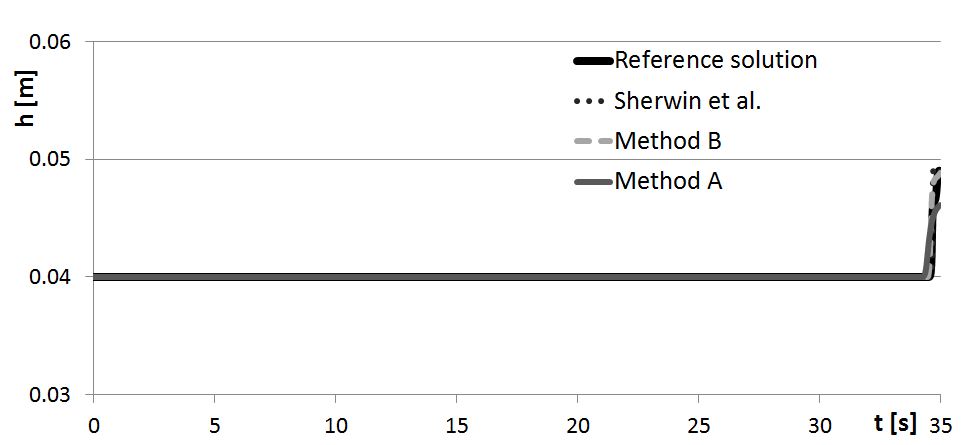}
\caption{Test 6 (subcritical wave): Computed free-surface elevation [m] in time [s] for Point 8.}
\label{rete_onda_grafici8}
\end{figure}
\begin{figure}[H]
\centering
\includegraphics[width=0.7\textwidth]{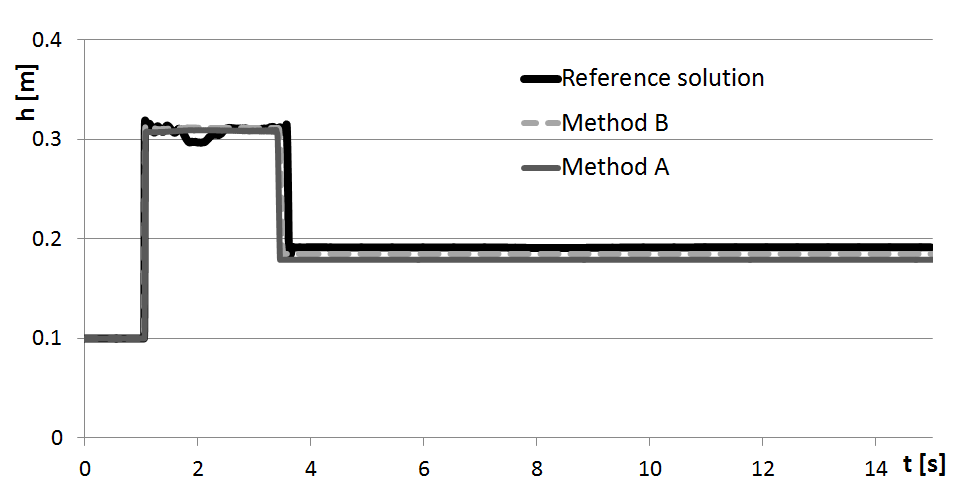}
\caption{Test 6 ( supercritical shock): Computed free-surface elevation [m] in time [s] for Point 1.}
\label{rete_shock_grafici1}
\end{figure}
\begin{figure}[H]
\centering
\includegraphics[width=0.7\textwidth]{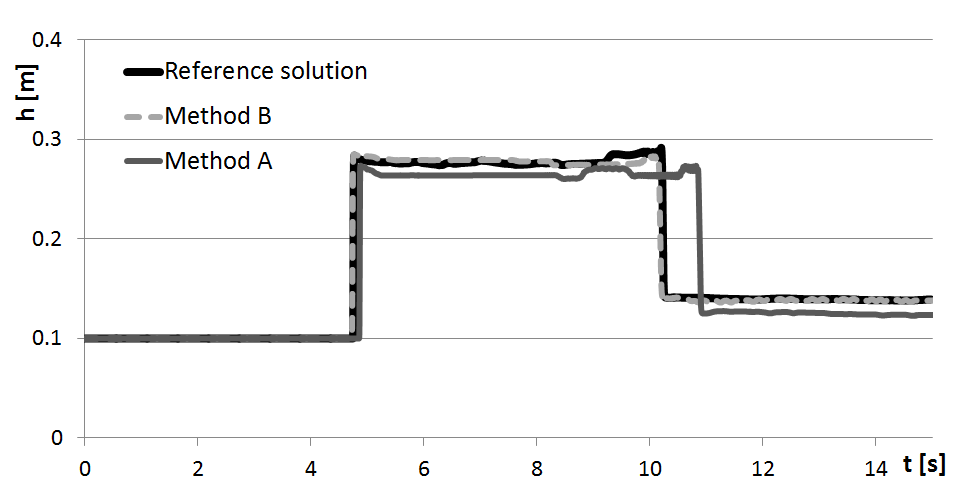}
\caption{Computed free-surface elevation [m] in time [s] for Point 3. Supercritical shock.}
\label{rete_shock_grafici3}
\end{figure}
\begin{figure}[H]
\centering
\includegraphics[width=0.7\textwidth]{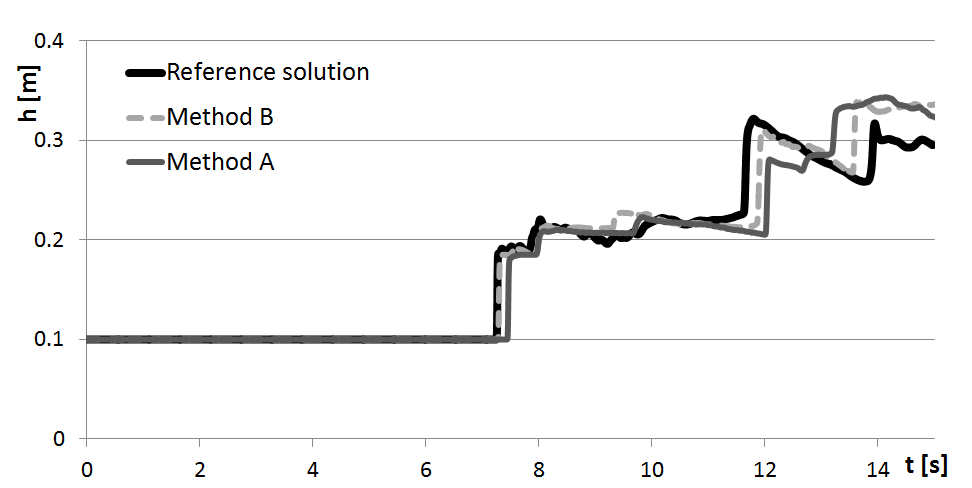}
\caption{Test 6 ( supercritical shock): Computed free-surface elevation [m] in time [s] for Point 6.}
\label{rete_shock_grafici6}
\end{figure}
\begin{figure}[H]
\centering
\includegraphics[width=0.7\textwidth]{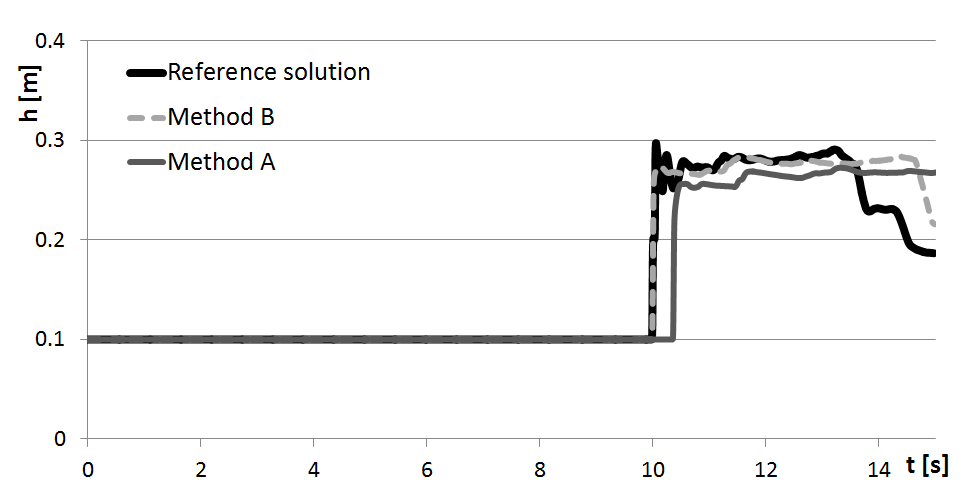}
\caption{Test 6 ( supercritical shock): Computed free-surface elevation [m] in time [s] for Point 8.}
\label{rete_shock_grafici8}
\end{figure}

\subsection{Computational times}

Here we show the computational times involved to solve each one of the six problems previously presented. Table \ref {tab_cputime}  shows the tests on the left columns and the CPU times in seconds in the subsequent columns, for the various methods used. Missing values for the PSFP method regard tests for which this method did not work. 
\begin{table}[H]
\renewcommand\arraystretch{1.2}
\centering {\small
\begin{tabular}{lcccc}
\hline
 & & & & \\[-4.5mm]
 {\bf Test} & {\bf 2D Reference} & {\bf PSFP method} & {\bf Method B} & {\bf Method A} \\
 & & & & \\[-4.5mm]
\hline
 & & & & \\[-4.5mm]
 Test 1  & 392.1 & 3.31 & 34.8 & 1.08 \\
 & & & & \\[-4.5mm]
Test 2 & 547.2 & 3.28 & 31.7 & 2.06 \\
 & & & & \\[-4.5mm]
Test 3 & 1215 & - & 68.9 & 3.53 \\
 & & & & \\[-4.5mm]
 Test 4 & 1684 & - & 128.3 & 4.69 \\
& & & & \\[-4.5mm]
 Test 6 (subcritical wave) & 5787 & 70.3 & 1413 & 19.2 \\
 & & & & \\[-4.5mm]
Test 6 (supercritical shock) & 13775 & - & 3091 & 51.5 \\
 & & & & \\[-4.5mm]
\hline
\end{tabular}}
\caption{Computational times [s] for all numerical methods reported in this paper, for six test problems.}
\label{tab_cputime}
\end{table}
As expected, the largest CPU times are those for the full 2D solver used to produce reference solutions. In terms of cost, there follows method B. The next one in CPU time  cost is the PSFP method, with Method A being the fastest, even faster that the simplest of all methods, namely the PSFP method, which is based entirely on 1D assumptions. It appears as if the method of choice is our method A, since it runs for all very demanding test problems, while giving reasonably accurate solutions as compared to the full 2D solver and at lowest computational cost. Computational saving factors for method A, relative to the full 2D solver, are of the order of 300, making the method a realistic option for complex applications.

\section{Concluding remarks} \label{sec:discussion} 

We have presented a novel method to treat junctions in networks of 1D shallow water channels. The method, called method A,  inserts a single 2D, junction-shaped  finite volume right at the junction, taking care that the element protrudes into  the 1D channels.  In this manner, the geometrical information, such as bifurcation angles and reflective boundaries is accounted for, locally. Method B results from generalising method A by inserting a local 2D unstructured grid in the vicinity of the junction. In addition, we briefly reviewed the existing junction method due to Peir\'o, Sherwin, Formaggia and Parker \cite{SherwinFormaggia,Sherwin2003}, which we termed PSFP method. All three approximate junction methods are assessed through a carefully selected suite of demanding test problems. No exact solutions to these problems exist to test the accuracy of approximate junction methods. We therefore use a fully 2D unstructured-mesh, second order method of the ADER type to compute accurate numerical solutions.  Method A is the preferred one, since it is simple and sufficiently accurate for all test problems. Method B is the more accurate of all three approximate methods tested, but also the most expensive, as shown by our computational efficiency test. Method A is the fastest, about three times faster than the  PSFP method and about 70 times faster than method B for the more realistic test problem involving a reasonably complex network. Methods A and B work well for all test problem, while the PSFP method only works for 3  of the 6 test problems. An attractive feature of method A, shared by method B, is that it can successfully cope with problems involving high subcritical, transcritical and supercritical flows at the junctions.  We note that due to the single-element of method A,  accuracy may deteriorate, depending on the mesh dimensions involved. This shortcoming is most evident in the first-order version of the methodology. Higher-order versions can ameliorate this deficiency. In fact, second order accuracy is found to be satisfactory, though we found a test problem for which only the third-order scheme produced fully satisfactory solutions, not shown here. Potential users of the schemes may have to assess this aspect of the methods before embarking on practical applications. For practical applications, both methods A and B may benefit from using local time-stepping, for example, following the methodology proposed in \cite{Dumbser:2007c,Mueller:2016}. This may be required by the disparity of spatial mesh sizes at the junctions and the 1D domains, which potentially implies disparity in time step sizes.

The methods presented in this paper can be applied to any problem involving networks of nearly straight 1D domains, provided the multidimensional version of the equations, 2D or 3D,  are available.

\vspace{10mm}

\begin{center}

{\bf Acknowledgements}

\end{center}

The authors are indebted to Prof. Dr. M. Dumbser, University of Trento, for useful discussions on the subject.

\newpage

\appendix
\section{Existing method for junctions: the PSFP method}\label{section_1Dexisting}

Here we briefly present an existing method for dealing with junctions, which we attribute to the work of J. Peir\'o, S. Sherwin, L. Formaggia and K. Parker, see  \cite{SherwinFormaggia} and \cite{Sherwin2003}. The methodology has become exceedingly popular and has been heavily applied in cardiovascular mathematics to deal with junctions of blood vessels. Here we adapt the method to deal with junctions in networks of water channels 

\subsection{Brief description of the method for shallow water}

\begin{figure}[H]
\centering
\subfigure[]{
\includegraphics[width=0.25\textwidth]{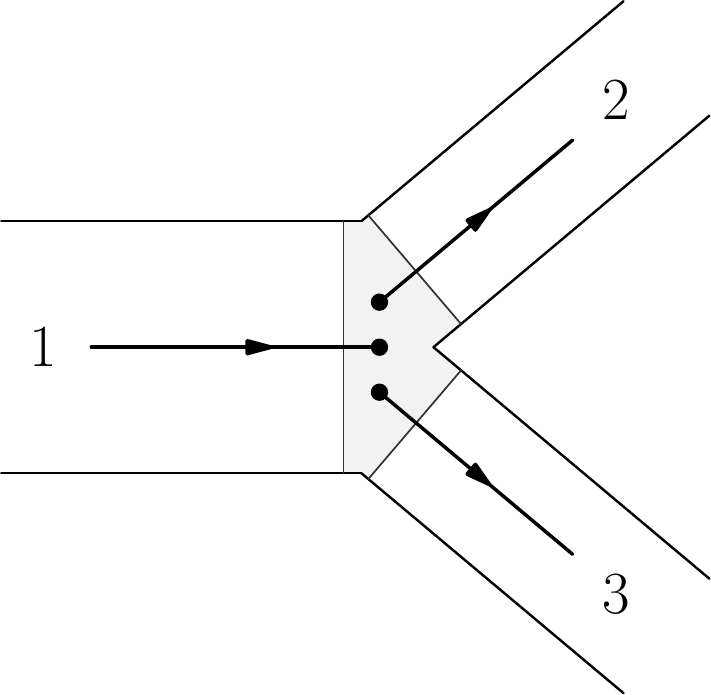}\label{schema_junc}}
\hspace{1cm}
\subfigure[]{
\includegraphics[width=0.29\textwidth]{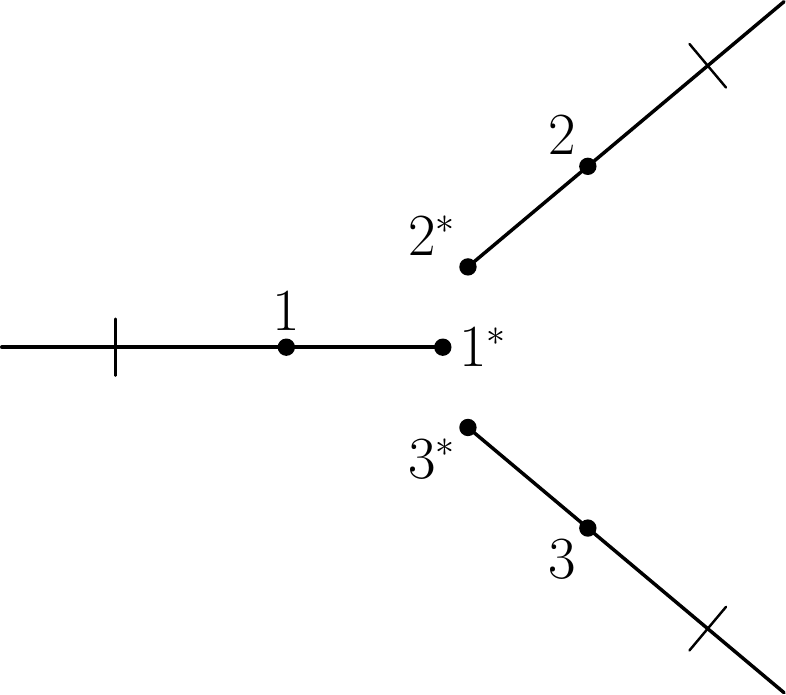}\label{schema_junc2}}
\caption{Junction problem as dealt with by the PSFP method. Frame (a): Idealised junction connecting a parent chanel 1 to two daughter channels 1 and 2. Frame (b): notation setup to solve the problem via 3 Riemann problems.}
\end{figure}

Consider the model bifurcation configuration shown in figure \ref{schema_junc}, where we denote the parent channel by index 1 and the daughter channels by  indices 2 and 3. We divide the domain composed by the three 1D channels into three sub-domains joined at the junction. In order to solve the junction problem and connect the 3 channels, a purely one-dimensional approach is taken to compute the needed 1D intercell numerical fluxes. 
At the bifurcation we have six unknowns, namely ($h_1^*$, $u_1^*$), ($h_2^*$, $u_2^*$) and ($h_3^*$, $u_3^*$). Once these are determined, corresponding numerical fluxes at the appropriate boundaries are computed, see figure \ref{schema_junc2}. Such fluxes will be used to update the 1D elements 1, 2 and 3.

The method assumes subcritical conditions and that the Riemann problem solutions involve only rarefaction waves. This allows the application of Riemann invariants, leading to 3 algebraic equations
\begin{equation}\label{prime3eq}
\begin{array}{l}
u_1^*+2\sqrt{g\,h_1^*}=u_1+2\sqrt{g\,h_1} \;,\\
\\
u_2^*-2\sqrt{g\,h_2^*}=u_2-2\sqrt{g\,h_2} \;, \\
\\
u_3^*\mp 2\sqrt{g\,h_3^*}=u_3\mp 2\sqrt{g\,h_3} \;,
\end{array}
\end{equation}
where the minus or plus sign in the third equation depends on the type of junction (diverging or merging).

Three more algebraic equations are
\begin{equation}\label{seconde3eq}
\begin{array}{l}
h_1^*u_1^*b_1=h_2^*u_2^*b_2+h_3^*u_3^*b_3 \;,\\
\\
h_1^*+\displaystyle{\frac{{u_1^*}^2}{2g}}=h_2^*+\displaystyle{\frac{{u_2^*}^2}{2g}} \;,\\
\\
h_1^*+\displaystyle{\frac{{u_1^*}^2}{2g}}=h_3^*+\displaystyle{\frac{{u_3^*}^2}{2g}} \;.
\end{array}
\end{equation}
The first of these represents  mass conservation  through the bifurcation, while the last two conditions are obtained from the requirement of continuity of total energy at the bifurcation. Conditions \eqref{prime3eq} and \eqref{seconde3eq} constitute a system of six non-linear  algebraic equations; these equations are  solved by a Newton method.

The described PSFP method requires the solution of a $ 6 \times 6$ non-linear algebraic system at each junction at each time step and represents a substantial computational effort in a complex network. However the method is indeed very simple and performs well under subcritical flow conditions with relatively small Froude numbers.  Sharing our experience, a word of caution is due. The method does not acknowledge the geometrical features of the problem and consequently, works better for small bifurcation angles and symmetric cases. For supercritical and subcritical flows with relatively high  Froude numbers the method may actually fail. Below we report a subcritical example with the following data:
\begin{equation}
\begin{array}{lll}
b_1=0.4\,m\;, & h_1=0.2\,m\;, & u_1=0.96\,m/s\;,\\
b_2=0.3\,m\;, & h_2=0.1\,m\;, & u_2=0.08\,m/s\;,\\
b_3=0.3\,m\;, & h_3=0.1\,m\;, & u_3=0.08\,m/s\;.
\end{array}
\end{equation}
Froude number in channel 1 is $Fr_1=0.685$ and in channel 2 and 3 is $Fr_2=Fr_3=0.081$, so the regime is clearly subcritical in every channel. Through the use of MAPLE we found two analytical solutions of system \eqref{prime3eq}-\eqref{seconde3eq}, namely:
\begin{equation}
\begin{array}{l}
h_1^*=0.155-0.0660I\,m\;,\\
h_2^*=0.186-0.0006I\,m\;,\\
h_3^*=0.186-0.0006I\,m\;,\\
u_1^*=1.240+0.5138I\,m/s\;,\\
u_2^*=0.807-0.0047I\,m/s\;,\\
u_3^*=0.807-0.0047I\,m/s\;,\\
\end{array}
\end{equation}
and
\begin{equation}
\begin{array}{l}
h_1^*=0.155+0.0660I\,m\;,\\
h_2^*=0.186+0.0006I\,m\;,\\
h_3^*=0.186+0.0006I\,m\;,\\
u_1^*=1.240-0.5138I\,m/s\;,\\
u_2^*=0.807+0.0047I\,m/s\;,\\
u_3^*=0.807+0.0047I\,m/s\;.\\
\end{array}
\end{equation}
Both solutions are complex;  needless to say, the Newton method did not converge.

\subsection{Test cases for the PSFP method}

In order to test the PSFP method we performed various tests consisting of a subcritical wave propagating across junctions of different geometry. We set initial conditions of zero velocity everywhere and water height of $16\,cm$. Boundary conditions at the outflow of daughter channels are transparent and velocity was assigned at the inflow of the parent channel by the relation
\begin{equation}
u(t)=0.4\exp(-0.5(t-3)^2)\;.
\end{equation}
To guarantee no wave reflection at the inflow, we used the method applied in  \cite{Alastruey}. Note that the flow regime of this test is subcritical, the maximum Froude number is about $0.5$.

\begin{figure}[H]
\centering
\addtocounter{subfigure}{-3}
\subfigure{
\includegraphics[width=0.33\textwidth]{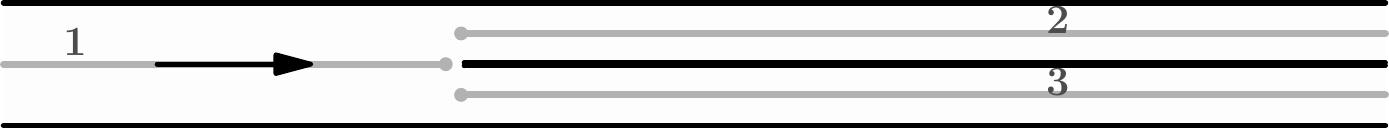}}
\subfigure{
\includegraphics[width=0.31\textwidth]{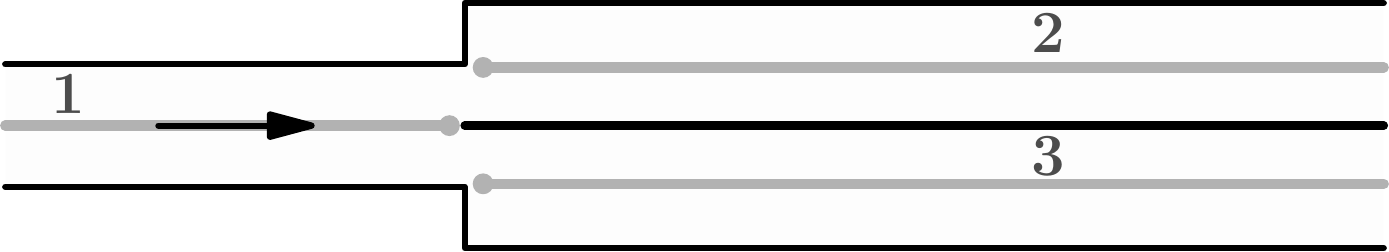}}
\subfigure{
\includegraphics[width=0.31\textwidth]{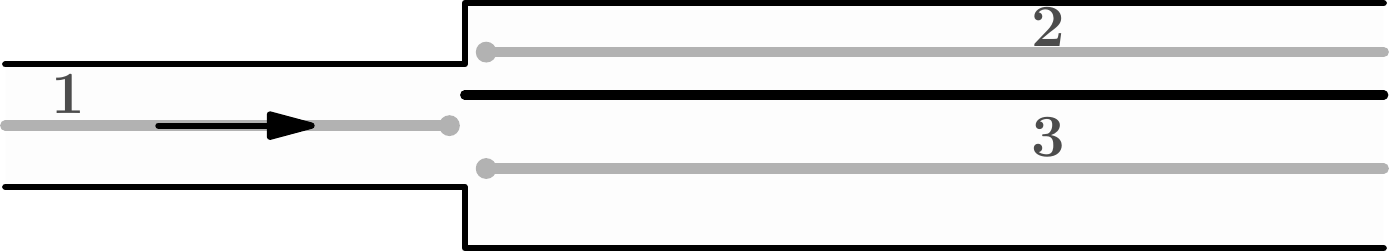}}
\subfigure[]{
\includegraphics[width=0.325\textwidth]{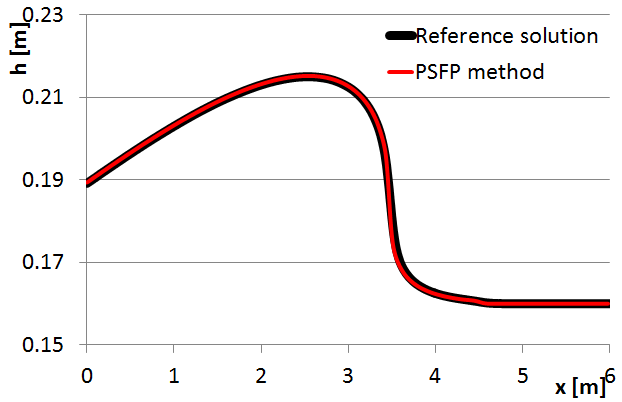}\label{canaledritto_t15}}
\subfigure[]{
\includegraphics[width=0.315\textwidth]{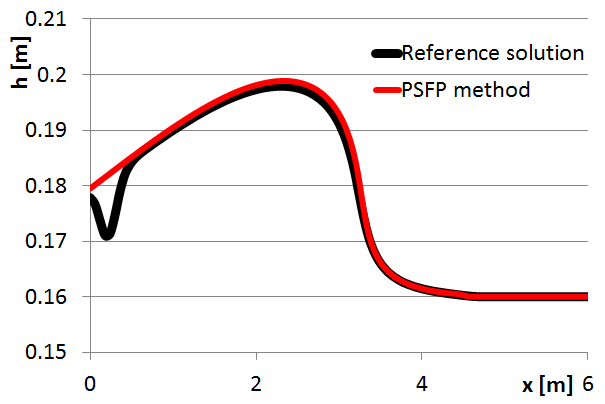}\label{canaledritto_allarg_2_t15}}
\subfigure[]{
\includegraphics[width=0.315\textwidth]{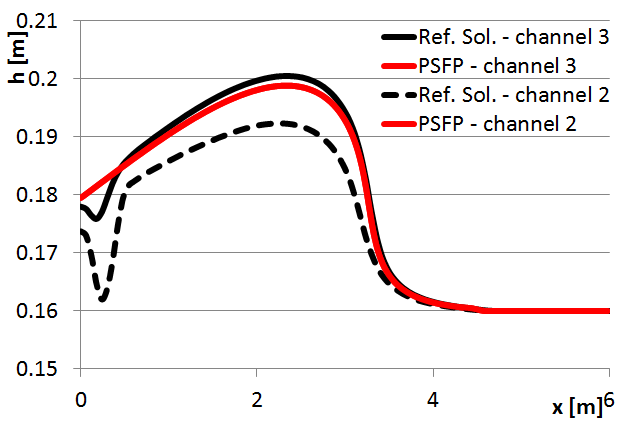}\label{canaledritto_allarg_as_t15}}
\caption{Performance of the PSFP junction method applied to problems of a subcritical wave propagating across straight channels. The PSFP-based 1D solutions are compared to the 2D reference solutions.}
\label{canaledritto}
\end{figure}

First we performed a test with a straight channel of width $0.4\,m$ which bifurcates into two channels of width $0.2\,m$ each placed at an angle of $0^\circ$, as dispayed in figure \ref{canaledritto_t15}. Then we tested the method with the same wave as above propagating through a channel which bifurcates into two larger channels placed at an angle of $0^\circ$, as displayed in figure \ref{canaledritto_allarg_2_t15}. Finally we performed a test with a situation similar to the previous one, but asymmetrical, as displayed in figure \ref{canaledritto_allarg_as_t15}. Figure \ref{canaledritto} shows computed results, comparing results with the reference 2D solution; water height profiles obtained in the secondary channels at time $t=6\,s$ are shown. Results worsen as the test problem depart from symetry.

\begin{figure}[H]
\centering
\includegraphics[width=0.4\textwidth]{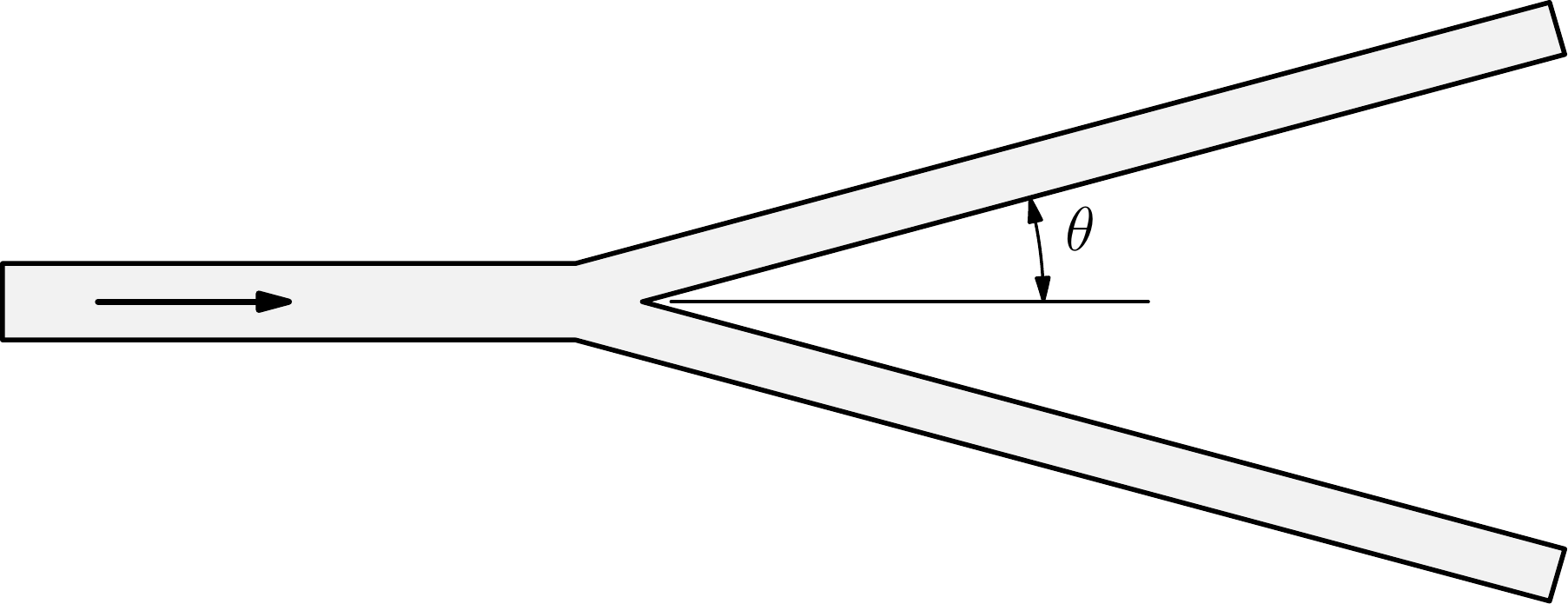}
\caption{Example of a symmetric bifurcation with bifurcation angle $\theta$.}
\label{junction15}
\end{figure}

We also tested the method for tests with bifurcation angles different from zero (figure \ref{junction15}). Initial conditions are the same as for the previous problems. 
We tested three bifurcations, with angles of $15^\circ$, $45^\circ$ and $90^\circ$. Figures \ref{canaledritto_t15}, \ref{canaledritto_allarg_2_t15} and \ref{canaledritto_allarg_as_t15} compare the water height profiles obtained in the primary channels at a time $t=8\,s$. Figures \ref{canaledritto_t15}, \ref{canaledritto_allarg_2_t15} and \ref{canaledritto_allarg_as_t15} compare the water height profiles in the secondary channels to the reference 2D solution.
\begin{figure}[H]
\centering
\subfigure[$\theta=15^\circ$, channel 1]{
\includegraphics[width=0.31\textwidth]{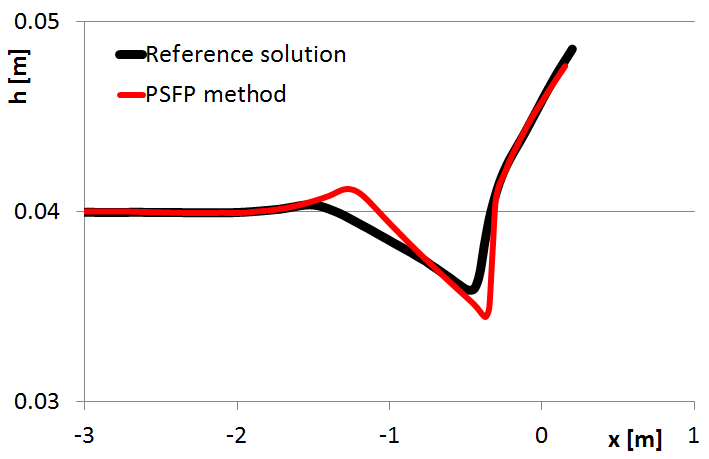}\label{junction15_2D1D_t20_ch1}}
\subfigure[$\theta=45^\circ$, channel 1]{
\includegraphics[width=0.31\textwidth]{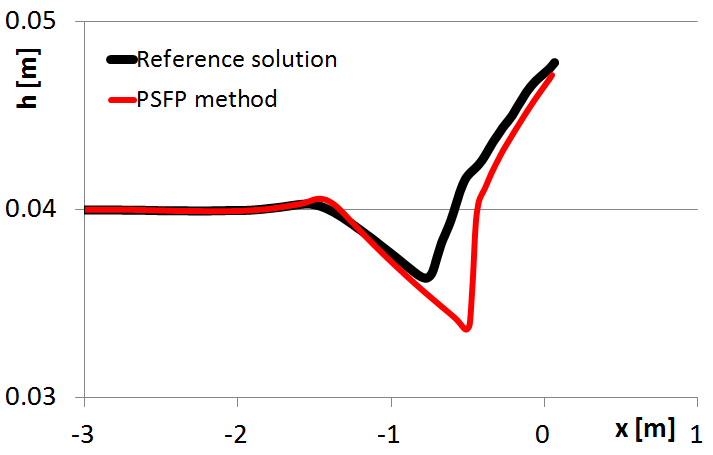}\label{junction45_2D1D_t20_ch1}}
\subfigure[$\theta=90^\circ$, channel 1]{
\includegraphics[width=0.31\textwidth]{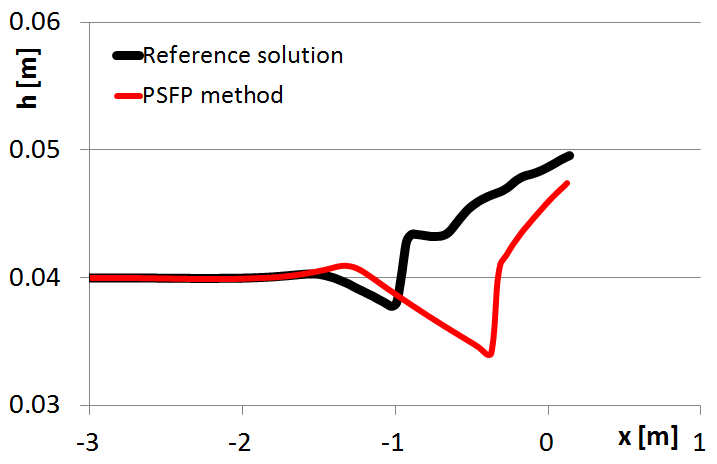}\label{junction90_2D1D_t20_ch1}}
\subfigure[$\theta=15^\circ$, channel 2]{
\includegraphics[width=0.31\textwidth]{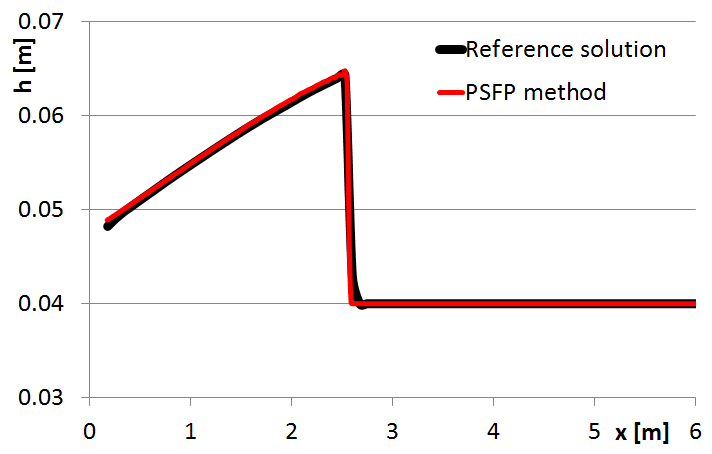}\label{junction15_2D1D_t20_ch2}}
\subfigure[$\theta=45^\circ$, channel 2]{
\includegraphics[width=0.31\textwidth]{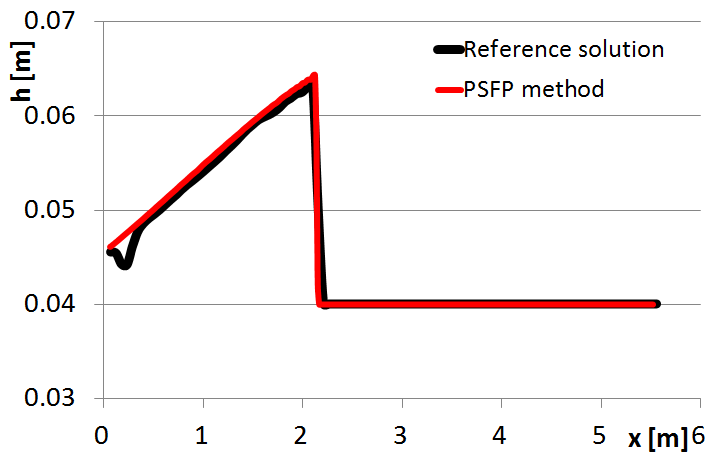}\label{junction45_2D1D_t20_ch2}}
\subfigure[$\theta=90^\circ$, channel 2]{
\includegraphics[width=0.31\textwidth]{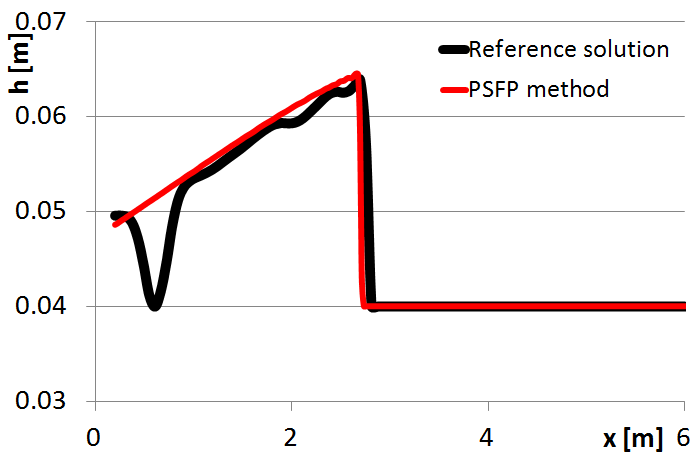}\label{junction90_2D1D_t20_ch2}}
\caption{Performance of the PSFP junction method applied to a problem of a  subcritical wave propagating across junctions of different bifurcation angles $\theta$. PSFP-based 1D solutions are compared to 2D reference solutions.}
\end{figure}
As the bifurcation angle increases, the solution in channel 1 obtained with the one-dimensional model based on the PSFP method, differs more and more from the reference two-dimensional solution. We note that solutions in channels 2 and 3 seem to be less sensitive to increased angle values.

\newpage

\section{A reference two-dimensional solver}\label{chapter_theoretical}

Here we give some details of our 2D shalllow water solver used to compute reference solutions to assess the junction models. Departing from equations (\ref{swe1}-\ref{swe2}) and making use of the rotational invariance  we construct a finite volume method on unstructured triangula meshes, written in the form

%
%
%
\begin{equation}                                 \label{fv}    
     {\bf Q}_{k}^{n+1} = {\bf Q}_{k}^{n} -\frac{\Delta t}{ |V_{k}| }\sum_{e=1}^{N}  {\cal{L}}_{e} {\bf T}_{e}^{-1} \hat{{\bf F}}_{e} \;.     		
\end{equation}
where the flux $\hat{{\bf F}}_{e}$ for edge $e$ is computed as discussed in Sect. 1.  The scheme is second-order accurate in space and time and its contruction follows the ADER approach \cite{ADER}. The first-order component of the method uses the HLLC flux \cite{HLLC}; the high-order version requires (i) a spatial nonlinear reconstruction procedure and (ii) the solution of the generalised Riemann problem. These steps are decribed in more detail below.\\

\noindent{\bf Spatial reconstruction.}
For the spatial reconstruction procedure on triangular meshes we follow \cite{AIAA89}; as a matter of fact we also use this procedure for 1D domains. For the 2D case, for the three triangles surrounding a general element we write the equation 
\begin{equation}
{\bf Q}_{J(k)}^n=a+b(x_{J(k)}-x_{j})+c(y_{J(k)}-y_{j})\;,
\end{equation}
where $j$ is the global index and $k$ the local index of each neighbour;  $x_{J(k)}$ and $y_{J(k)}$ the coordinates of the centroids of the neighbours. Once this three-equation system is solved, in order to respect the conservation property we have to force the coefficient $a$ to satisfy $a={\bf Q}_{i}^n$. Coefficients $b$ and $c$ represent the best estimates of the solution gradient. We have to impose a limiter on this gradient in order to respect the non-oscillatory property. Barth and Jespersen \cite{AIAA89} suggest to find the largest admissible $\Phi$ so that the values of the linearly reconstructed function do not exceed the maximum and minimum of neighbouring centroid values (including the centroid value in the element). Then, keeping in mind that for linear reconstructions extrema in a triangle occur at the vertices, we can calculate the limiter $\Phi$ as the minimum $\phi_k$ found among all vertices, and finally multiply $b$ and $c$ by $\Phi$. We perform this procedure for every conserved variable. \\

\noindent{\bf Generalised Riemann problem.}
There exist several methods to solve the Generalised Riemann Problem (GRP), some of which are described in the paper by Castro and Toro \cite{CastroToro2008}. Here we have implemented the method proposed by Harten and collaboartors in \cite{Harten}, reinterpreted in \cite{CastroToro2008} as a solver for the GRP. The method evolves the Riemann problem data ${\bf Q}_L$ and ${\bf Q}_R$ considering power series expansions in time on each side of the interface. For a second order scheme we need first order time series expansions, as follows:
\begin{equation}
\begin{array}{l}
\widetilde{{\bf Q}}_L(\frac{\Delta t}{2})={\bf Q}_L(0_-)+\frac{\Delta t}{2}\partial_t{\bf Q}(0_-,0)\;,\\
\\
\widetilde{{\bf Q}}_R(\frac{\Delta t}{2})={\bf Q}_R(0_+)+\frac{\Delta t}{2}\partial_t{\bf Q}(0_+,0)\;.
\end{array}
\end{equation}
Using the Cauchy-Kowalewskaya procedure to express time derivatives as functions of space derivatives, and given the slopes resulting from second order reconstruction, we have:
\begin{equation}
\begin{array}{l}
\partial_{t}{\bf Q}(0_-,0) = - {\bf A}({\bf Q}(0_-,0))\cdot {\bf b}({\bf Q}(0_-,0)) - {\bf B}({\bf Q}(0_-,0))\cdot {\bf c}({\bf Q}(0_-,0)) \;,\\
\\
\partial_{t}{\bf Q}(0_+,0) = - {\bf A}({\bf Q}(0_+,0))\cdot {\bf b}({\bf Q}(0_+,0)) - {\bf B}({\bf Q}(0_+,0))\cdot {\bf c}({\bf Q}(0_+,0)) \;.
\end{array}
\end{equation}
Then we apply the rotational matrix ${\bf T}_s$ on $\widetilde{{\bf Q}}_L$ and $\widetilde{{\bf Q}}_R$, obtaining  $\hat{\widetilde{{\bf Q}}}_{L}$ and $\hat{\widetilde{{\bf Q}}}_{R}$, and we compute the numerical flux from the solution of a classical Riemann problem, for which we use the HLLC approximation \cite{HLLC}. The resulting method is second-order accurate in space and time.\\

\noindent{\bf Grid-independence test.}
In order to choose a mesh size that gives a reliable reference solution, we performed a 
grid-independence test. We used a typical problem of a shock wave propagating across a junction. For this problem we stored the computed free-surface elevation in time for a representative point placed where the two-dimensional character of the problem is more evident. Then we run the simulation with five different mesh sizes, each time halving the element size.  We calculated the integral in time for each mesh and then we computed the change relative to the previous mesh. Results are shown in table \ref{Grid_independence_tab}. We considere the difference obtained for the mesh of size $0.02$  as beeing satisfactory. For all the computations to provide reference 2D solutions we used this mesh size.
\begin{table}[H]
\renewcommand\arraystretch{1.2}
\centering{\small
\begin{tabular}{cccc}
\hline
 & & \\[-4.5mm]
{\bf Grid} & {\bf Integrals} & {\bf Relative differences} & {\bf CPU times} \\
 & & \\[-4.5mm]
\hline
 & & \\[-4.5mm]
 0.16 & 0.3767 & 0.0502 & 2.8 s \\
 & & \\[-4.5mm]
 0.08 & 0.3966 & 0.1216 & 22.3 s \\
 & & \\[-4.5mm]
 0.04 & 0.4515 & 0.0046 & 129.1 s \\
 & & \\[-4.5mm]
 0.02 & 0.4536 & 0.0015 & 1215 s \\
 & & \\[-4.5mm]
 0.01 & 0.4543 &  & 9528 s \\
 & & \\[-4.5mm]
\hline
\end{tabular}  }
\caption{Variation of integrals of water height as the mesh is refined.}
\label{Grid_independence_tab}
\end{table}

\cleardoublepage
\addcontentsline{toc}{chapter}{\bibname}

\bibliographystyle{myplainnat}
\bibliography{bibliografia}

\end{document}